\documentclass[a4paper,12pt]{article}
\usepackage[text={6.5in,9in},centering]{geometry}
\usepackage{graphicx,url,subfig}
\usepackage{multirow}
\RequirePackage[OT1]{fontenc}
\RequirePackage{amsthm,amsmath,amssymb,amscd}
\RequirePackage[numbers]{natbib}
\RequirePackage[colorlinks,citecolor=blue,urlcolor
=blue]{hyperref}
\usepackage{enumerate}
\usepackage{datetime}
\usepackage{etex}
\makeatother
\numberwithin{equation}{section}
\allowdisplaybreaks
\usepackage{accents}
\newcommand{\ubar}[1]{\underaccent{\bar}{#1}}

\theoremstyle{plain}
\newtheorem{theorem}{Theorem}[section]

\newtheorem{lemma}[theorem]{Lemma}
\newtheorem{proposition}[theorem]{Proposition}

\theoremstyle{definition}
\newtheorem{definition}[theorem]{Definition}

\newtheorem{assumption}[theorem]{Assumption}

\newtheorem{remark}[theorem]{Remark}

\newtheorem{example}[theorem]{Example}

\newcommand{\E}{\mathbb{E}}
\newcommand{\W}{\dot{W}}
\newcommand{\ud}{\ensuremath{\mathrm{d}}}
\newcommand{\Ceil}[1]{\left\lceil #1 \right\rceil}

\newcommand{\Indt}[1]{1_{\left\{#1 \right\}}}

\newcommand{\Norm}[1]{\left|\left|  #1   \right|\right|}

\newcommand{\Id}{\text{Id}}

\newcommand{\calB}{\mathcal{B}}

\newcommand{\calF}{\mathcal{F}}
\newcommand{\calG}{\mathcal{G}}
\newcommand{\calK}{\mathcal{K}}

\newcommand{\calL}{\mathcal{L}}
\newcommand{\calM}{\mathcal{M}}
\newcommand{\calN}{\mathcal{N}}

\newcommand{\G}{\mathcal{G}}
\newcommand{\calP}{\mathcal{P}}

\newcommand{\bbC}{\mathbb{C}}
\newcommand{\bbN}{\mathbb{N}}

\newcommand{\R}{\mathbb{R}}
\newcommand{\myEnd}{\hfill$\square$}

\newcommand{\Erf}{\ensuremath{\mathrm{erf}}}
\newcommand{\Erfc}{\ensuremath{\mathrm{erfc}}}

\DeclareMathOperator{\Lip}{\mathit{L}}
\DeclareMathOperator{\LIP}{Lip}
\DeclareMathOperator{\lip}{\mathit{l}}

\DeclareMathOperator{\Vip}{\overline{\varsigma}}
\DeclareMathOperator{\vip}{\underline{\varsigma}}

\newcommand{\Dxa}{{}_x D_\delta^a}
\newcommand{\DCap}{{}_t D_*^\beta}
\newcommand{\DRL}[2]{{}_t D_{#1}^{#2}}
\newcommand{\lMr}[3]{\:{}_{#1} #2_{#3}}

\newcommand{\myRef}[2]{#1}

\title{Nonlinear stochastic time-fractional diffusion equations on $\R$:
moments, H\"older regularity and intermittency}

\author{\bf Le Chen\footnote{Research supported by both University of Utah and the Swiss National Foundation for Scientific
Research through an SNSF fellowship.} \footnote{Address:
Department of Mathematics, University of Utah, 155 S 1400 E RM 233, Salt Lake City, UT, 84112-0090,
USA. \textit{Email:} \url{chen@math.utah.edu}}\\[0.5em]
\it University of Utah
\date{}
}

\begin{document}
\maketitle

\begin{center}
\begin{minipage}[rct]{5 in}
\footnotesize \textbf{Abstract:}
We study the nonlinear stochastic time-fractional
diffusion equations in the spatial domain $\R$, driven by multiplicative space-time
white noise.  The fractional index
$\beta$ varies continuously from $0$ to $2$.
The case $\beta=1$ (resp. $\beta=2$) corresponds
to the stochastic heat (resp.
wave) equation. The cases $\beta\in \:]0,1[\:$
and $\beta\in \:]1,2[\:$ are
called {\it slow diffusion equations} and {\it fast diffusion equations}, respectively.
Existence and uniqueness of random field solutions with measure-valued initial data, such as the Dirac delta measure, are established.
Upper bounds on all $p$-th moments $(p\ge 2)$ are obtained, which
are expressed using a kernel function $\calK(t,x)$.
The second moment is sharp.
We obtain the H\"older continuity of the solution for the slow diffusion equations
when the initial data is a bounded function.
We prove the weak intermittency for both slow and fast diffusion equations.
In this study, we introduce a special function, the {\it two-parameter Mainardi functions}, which are generalizations of the one-parameter Mainardi functions.

\vspace{2ex}

\vspace{2ex}
\vspace{4ex}
\end{minipage}
\end{center}


\section{Introduction}
{\it Viscoelasticity} is
the property of materials that exhibit both viscous
and elastic characteristics when undergoing deformation (see e.g.
\cite{lemaitre1994mechanics,dimitrienko2010nonlinear,perkins2011viscoelasticity}
). Viscosity mainly refers to fluids and elasticity to solids.
A {\it linear} theory to bring these two properties together has been achieved using fractional calculus by Mainardi and his coauthors; see \cite{Mainardi10Book} for an introduction to this subject.
It has wide applications to fields such as chemistry (e.g. \cite{doi1988theory,ferry1980viscoelastic}), seismology (e.g. \cite{aki2002quantitative}), soil mechanics (e.g.
\cite{klausner1991fundamentals}), arterial rheology
(\cite{craiem2008fractional}), biological tissues (e.g. \cite{magin2010fractional}), etc.
In this linear model, the system is governed by the partial differential operator
\[
\calL=\lMr{x}{D}{\delta}^a - \lMr{t}{D}{*}^\beta,
\]
where the space derivative $\Dxa$ is the Riesz-Feller fractional
derivative of order $a$ and
skewness $\delta$, and the time
derivative $\lMr{t}{D}{*}^\beta$
is a {\it Caputo derivative} of order $\beta\in \: ]0,2]$.
These three parameters vary in the following ranges:
\[
a\in \:]0,1]\:,\quad \beta \in \:]0,2]\:,\quad |\delta| \le a\wedge (2-a)\wedge (2-\beta)\:,
\]
where $a\wedge b:=\min(a,b)$.
We are interested in this linear model driven by multiplicative
space-time white noise:
\[
\calL u(t,x) = I_t^{\Ceil{\beta}-\beta} \left[\rho(u(t,x)) \W(t,x)\right],\quad
t\in \R_+^*:=\;]0,+\infty[\;,\:x\in\R,
\]
where $\Ceil{\beta}$ be the smallest integer not less than $\beta$, $\dot{W}$ is the space-time white noise,
the function $\rho:\R\mapsto\R$ is Lipschitz continuous,
and $I_t^\alpha$ is the {\it Riemann-Liouville fractional integral} of order $\alpha$:
\[
I_t^\alpha f(t):= \frac{1}{\Gamma(\alpha)}\int_0^t (t-s)^{\alpha-1}f(s)\ud s,  \quad \text{for $t>0$ and $\alpha\ge 0$}.
\]
Let $\Id$ denote the identity operator.
When $\beta$ is an integer, then $ I_t^{\Ceil{\beta}-\beta}= I_t^0=\Id$.
The case where $\rho(u)=\lambda u$, $\beta=1$, and $a=2$ (hence $\delta=0$), is called
the {\it parabolic Anderson model}; see \cite{BertiniCancrini94Intermittence,CarmonaMolchanov94}.
The logarithm of the solution gives the {\em Hopf-Cole solution} to the famous
{\em Kardar-Parisi-Zhang equation} \cite{KPZ86}.

Due to the time-fractional derivative, the semigroup theory does not work except for the case $\beta=1$.
These studies heavily depend on the properties of the fundamental solutions or
the Green functions to $\calL u =0$.
In \cite{MainardiEtc01Fundamental}, some of these
Green functions are
obtained through inverse Fourier transform of some special functions, among which the
following three cases are
more trackable:
(1) Space-fractional heat equation: $\{0<a\le 2,\: \beta =1\}$;
(2) Time-fractional heat/wave equation: $\{a= 2,\: 0< \beta \le 2\}$;
(3) Neutral fractional diffusion equation: $\{0<a= \beta \le 2\}$.
The first case has been recently studied in
\cite{ChenDalang14FracHeat,ChenKim14Comparison}.
In this paper, we will study the second case, i.e., we will study the following nonlinear
stochastic time-fractional diffusion equations (formally):
\begin{align}\label{E:tFracHt}
  \left(\displaystyle \DCap -
\frac{\partial^2}{\partial x^2} \right) u(t,x) =
I_t^{\Ceil{\beta}-\beta} \left[\rho\left(u(t,x)\right) \dot{W}(t,x)\right], \quad \beta \in \;]0,2],\: t\in \R_+^*,\:
x\in\R.
\end{align}
The Caputo fractional differential operator $\DCap$ is defined as
\[
\DCap f(t) :=
\begin{cases}
\displaystyle
 \frac{1}{\Gamma(m-\beta)} \int_0^t\ud
\tau\: \frac{f^{(m)}(\tau)}{(t-\tau)^{\beta+1-m}}&
\text{if $m-1<\beta<m$\;,}\\[1em]
\displaystyle
\frac{\ud^m}{\ud t^m}f(t)& \text{if $\beta=m$}\;.
\end{cases}
\]
We refer to \cite{Diethelm10AFDE,MainardiEtc01Fundamental} for more
details of these fractional differential operators.
When $\beta=2$,
$\DCap=\frac{\partial^2}{\partial t^2}$ and
\eqref{E:tFracHt} reduces to the
{\it stochastic wave equation} (SWE):
\begin{align}
 \label{E:Wave}
\left(\frac{\partial^2}{\partial t^2}- \kappa^2
\frac{\partial^2}{\partial x^2}\right) u(t,x) = \rho(u(t,x)) \W(t,x)\:,
\end{align}
with the speed of wave propagation $\kappa=1$.
When $\beta=1$, $\DCap=\frac{\partial}{\partial t}$
and \eqref{E:tFracHt} reduces to the
{\it stochastic heat equation} (SHE):
\begin{align}
 \label{E:Heat}
\left(\frac{\partial}{\partial t}- \frac{\nu}{2}
\frac{\partial^2}{\partial x^2}\right) u(t,x) = \rho(u(t,x)) \W(t,x)\:,
\end{align}
with the diffusion parameter $\nu=2$.
The above two special cases have been studied carefully; see \cite{LeChen13Thesis, ChenDalang13Heat,ChenDalang13Holder,ChenDalang14Wave,ConusEct12Initial}.
The case $\beta\in \;]0,1]$ is called the {\it slow diffusion},
$\beta\in\;]1,2]$ the {\it fast diffusion}, and $\beta=1$
the {\it standard diffusion}. In the following
we will also call the case $\beta\in\:]0,1[\:$ slow diffusion
and the case $\beta\in\:]1,2[\:$ fast diffusion.
For the slow and standard diffusions, we only need to specify the initial data $u(0,x)$.
For the fast diffusion, we need to give $\frac{\partial}{\partial t}u(0,x)$ as well.
Note that another related equation is the {\it stochastic fractional heat equation} (SFHE):
\begin{align}
 \label{E:FracHeat}
\left(\frac{\partial}{\partial t}- \lMr{x}{D}{\delta}^a\right) u(t,x) = \rho(u(t,x)) \W(t,x)\:,
\end{align}
which has been studied recently in \cite{ChenDalang14FracHeat,ChenKim14Comparison};
see also \cite{DebbiDozzi05On,FoondunKhoshnevisan08Intermittence}.

All investigations on SPDEs of the above kinds  require a good study of the corresponding Green functions.
By Green functions, we mean the solutions to the following equations
\begin{align}\label{E:Homo-Eq}
 \begin{cases}
  \left(\displaystyle \DCap -
\frac{\partial^2}{\partial x^2} \right) u(t,x) =
0,& t\in \R_+^*,\:
x\in\R\:,
\\[1em]
\displaystyle
u(0,x) = \delta_0(x)\;,& x\in\R\;,\quad\text{if $0<\beta\le 1$\;,}\\[0.4em]
\displaystyle
u(0,x) = 0\;,\quad \frac{\partial }{\partial t}
u(0,x) = \delta_0(x)\;,&x\in\R\;,\quad \text{if $1<\beta\le 2$\;,}
 \end{cases}
\end{align}
where $\delta_0$ is the Dirac delta function with
a unit mass at zero. We use $G_\beta(t,x)$ to denote these Green functions.
The Green functions for slow diffusion equations and their properties can be found in \cite{MainardiEtc01Fundamental}.
As far as we know, there is no literature studying the Green functions of the fast diffusion equations.
Note  that in \cite{MainardiEtc01Fundamental}, the Green functions for the fast diffusions are not the one we need.
To obtain the Green functions for the fast diffusion equations, one needs to generalize the one-parameter {\it Mainardi function} (see \cite{Mainardi10Book,MainardiEtc01Fundamental}) to the
two-parameter settings (see \eqref{E:2p-Mainardi} below),
based on which corresponding properties for the Green functions of the fast diffusion equations need to be proved (see Lemma \ref{L:Green} below).

If we denote the solution to the homogeneous equation by $J_0(t,x)$ (see \eqref{E:J0} below),
then the rigorous meaning of \eqref{E:tFracHt},
which is the actual equation that we are going to study, is the following stochastic integral equation:
\begin{equation}
\label{E:WalshSI}
 \begin{aligned}
  u(t,x) &= J_0(t,x)+I(t,x),\quad\text{where}\cr
I(t,x) &=\iint_{[0,t]\times\R}
G_\beta\left(t-s,x-y\right)\rho\left(u(s,y)\right)
W(\ud s,\ud y),
 \end{aligned}
\end{equation}
where the stochastic integral is the Walsh integral \cite{Walsh86}.
To motivate the relation between the SPDE \eqref{E:tFracHt} and the integral equation \eqref{E:WalshSI},
we need the {\it time-fractional Duhamel's principle} (see \cite[Theorem 3.6]{UmarovSaydamatov06}).
If we replace the right hand side of \eqref{E:tFracHt} by a nice forcing term $g(t,x)$.
Then by \cite[Theorem 3.6]{UmarovSaydamatov06}, the solution to \eqref{E:tFracHt} 
with vanishing initial conditions
$\frac{\partial^m }{\partial t^m} u(0,x)=0$ for $m=0, \dots, \Ceil{\beta}-1$ is
\[
u(t,x)=\int_0^t\ud s\int_\R \ud y\: 
G_\beta(t-s,x-y) \: \DRL{+}{\Ceil{\beta}-\beta} g(s,y),
\]
where $\DRL{+}{\alpha}$ for $\alpha\ge 0$ is the {\it Riemann-Liouville fractional derivatives} of order $\alpha$:
\[
\DRL{+}{\alpha} f(t):=
\begin{cases}
\displaystyle
 \frac{1}{\Gamma(m-\alpha)} \frac{\ud^m}{\ud t^m}\int_0^t\ud
\tau\: \frac{f(\tau)}{(t-\tau)^{\alpha+1-m}}&
\text{if $m-1<\alpha<m$\;,}\\[1em]
\displaystyle
\frac{\ud^m}{\ud t^m}f(t)& \text{if $\alpha=m$}\;.
\end{cases}
\]
Then if one replaces $g(t,x)$ by $I_t^{\Ceil{\beta}-\beta} g(t,x)$ and uses the fact that $\DRL{+}{\alpha}\circ I_t^{\alpha}= \Id$ for all $\alpha\ge 0$, then one can see that solution to
\[
\begin{cases}
   \left(\displaystyle \DCap -
\frac{\partial^2}{\partial x^2} \right) u(t,x) =
I_t^{\Ceil{\beta}-\beta} \left[g(t,x)\right], & \beta \in \;]0,2],\: t\in \R_+^*,\:
x\in\R,\\[1em]
\displaystyle \frac{\partial^m}{\partial t^m} u(0,x) =0,& \text{for $m=0,\dots,\Ceil{\beta}-1$},
\end{cases}
\]
is
\[
u(t,x)=\int_0^t\ud s\int_\R \ud y\:
G_\beta(t-s,x-y) \: g(s,y).
\]

We will establish the existence of random field solutions to \eqref{E:WalshSI} starting from measure-valued initial conditions. Let $\mu$ be a Borel measure and $\mu=\mu_+-\mu_-$, where,
from the Jordan decomposition,  $\mu_\pm$ are two
nonnegative Borel measures with disjoint support and $|\mu|=\mu_++\mu_-$.
Define an axillary function
\begin{align}\label{E:f}
f_\eta(x) :=
\exp\left(-\frac{\eta}{2}|x|^{\frac{2}{2-\beta}}\right)\;,\quad\text{for $x\in\R$.}
\end{align}
Let $\calM(\R)$ be the set of signed (regular) Borel measures on $\R$.
For $0<\beta<2$, define
\begin{align}
 \label{E:InitData}
\calM_T^\beta \left(\R\right):=\Big\{
	\mu\in\calM(\R):
\left(|\mu|*f_\eta\right)(x)<+\infty\:,
\text{for all $\eta>0$ and $x\in\R$}
\Big\},
\end{align}
where $*$ denotes the convolution in the space variable.
Then $\calM_T^1(\R) = \calM_H(\R)$, where $\calM_H(\R)$
is the notation used in
\cite{ChenDalang13Heat, ChenDalang13Holder} for the admissible initial data for the SHE \eqref{E:Heat}.
Note that even though the initial data can be Schwartz distributions for the heat equation without noise, but for the SPDE, initial data cannot go beyond measures; see \cite[Theorem
3.2.17]{LeChen13Thesis} or \cite[Theorem 2.22]{ChenDalang14Wave}.
We will prove the existence and uniqueness of random field solutions to
\eqref{E:tFracHt} for all initial data in $\calM_T^\beta(\R)$.
As in \cite{ChenDalang13Heat,ChenDalang14Wave,ChenDalang14FracHeat},
we will obtain similar moment formulas expressed using a special
function $\calK(t,x)$. For the SHE and the SWE,
this kernel function $\calK(t,x)$ has an explicit form.
But for the space-fractional heat equations
\cite{ChenDalang14FracHeat} and the current time-fractional diffusion equations,
we only have some estimates on it.
In particular,  for the slow diffusion equations,
we will obtain both upper and lower bounds on $\calK(t,x)$. For the fast diffusion equations,
we will only derive some upper bounds.

After establishing the existence of random field solutions, we will study some properties of the solutions.
The first property is the sample-path regularity (for the slow diffusion equations).
Given a subset $D\subseteq \R_+\times\R$ and positive constants
$\beta_1,\beta_2$, denote by $C_{\beta_1,\beta_2}(D)$ the set of functions
$v: \R_+ \times \R \to \R$ with the following property:
for each compact set $K\subseteq D$,
there is a finite constant $C$ such that for all $(t,x)$ and $(s,y)\in K$,
\[
\vert v(t,x) - v(s,y) \vert \leq C \left(\vert
t-s\vert^{\beta_1} + \vert x-y \vert^{\beta_2}\right).
\]
Denote
\[
C_{\beta_1-,\beta_2-}(D) := \cap_{\alpha_1\in \;\left]0,\beta_1\right[}
\cap_{\alpha_2\in \;\left]0,\beta_2\right[} C_{\alpha_1,\alpha_2}(D)\;.
\]
We will show that for the slow diffusion equations,
if the initial data is a bounded function, i.e., $\mu(\ud x)=f(x)\ud x$ with
$f\in L^\infty(\R)$, then
\begin{align}\label{E:Holder-u}
u(\cdot,\circ) \in
C_{\frac{2-\beta}{4}-,\:\frac{1}{2}-}\left(\R_+\times\R\right),\quad \text{a.s.}
\end{align}
Moreover, if $f$ is bounded and
$\alpha$-H\"older continuous ($\alpha\in \;]0,1[\:$), then
\begin{align}\label{E:HHolder-u}
u(\cdot,\circ)\in
C_{\frac{2-\beta}{4}-,\frac{1}{2}-}\left(\R_+^*\times\R\right)
\;\cap\;
 C_{\left(\frac{ \alpha \beta}{2}\wedge\frac{2-\beta}{4}\right)-,\:
\left(\alpha\beta \wedge \frac{1}{2}\right)-}\left(\R_+\times\R\right)
\;,\quad\text{a.s.}
\end{align}
When $\beta=1$, the above results partially recover the results for the
stochastic heat equation in \cite{ChenDalang13Holder}.
Note that the regularity results in \cite{ChenDalang13Holder} is
more general since the initial data can be measures.

The second property that we are going to study is the intermittency.
More precisely, define the {\it upper and lower (moment) Lyapunov exponents} as follows
\begin{align}
\label{E:Lyapunov}
 \overline{m}_p(x):=& \mathop{\lim\sup}_{t\rightarrow+\infty} \;
\frac{\log\E\left[|u(t,x)|^p\right]}{t}\quad\text{and}\quad
\underline{m}_p(x):=\mathop{\lim\inf}_{t\rightarrow+\infty} \;
\frac{\log\E\left[|u(t,x)|^p\right]}{t}.
\end{align}
When the initial data are spatially homogeneous (i.e., the initial data are constants),
so is the solution $u(t,x)$,
and then the Lyapunov exponents do not depend on the spatial variable.
In this case, a solution is called {\it fully intermittent} if $m_1=0$ and
$\underline{m}_2>0$ (see \cite[Definition III.1.1, on p. 55]{CarmonaMolchanov94}).
As for the weak intermittency, there are various definitions.
For convenience of stating our results, we will call the solution {\it weakly intermittent of type I} if $\underline{m}_2>0$,
and {\it weakly intermittent of type II} if $\overline{m}_2>0$.
Clearly, the weak intermittency of type I is slightly stronger than the
the weak intermittency of type II, but weaker than the full intermittency by missing $m_1=0$.
The weak intermittency of type II is used in \cite{FoondunKhoshnevisan08Intermittence}.

The full intermittency for the SHE and the SFHE,
and the weak intermittency of type I for SWE are established in
\cite{BertiniCancrini94Intermittence}, \cite{ChenKim14Comparison} and
\cite{ChenDalang14Wave}, respectively.
Conus, {\it et al. } prove the weak intermittency of type II for the SWE in \cite[Theorem 2.3]{ConusEtal13Wave}.
We will establish the weak intermittency of type I for the slow diffusion equations and
the weak intermittency of type II for the fast diffusion equations.
Moreover, we show that
\begin{align}\label{E:WI}
 \overline{m}_p\le
\begin{cases}
 C\: p^{\frac{4-\beta}{2-\beta}} & \text{if $\beta\in \:]0,1]$,}\cr
C\: p^{\frac{8-\beta}{6-\beta}} & \text{if $\beta\in \:]1,2[\:$,}\cr
\end{cases}
\end{align}
which reduces to the SHE case (see
\cite{BertiniCancrini94Intermittence,ChenDalang13Heat,FoondunKhoshnevisan08Intermittence})
when $\beta =1$, i.e., $\overline{m}_p \le C\: p^3$,
and to the SWE case (see \cite{ChenDalang14Wave}) when $\beta=2$, i.e.,
$\overline{m}_p \le C\:p^{3/2}$.
Note that the above constants $C$ may vary from one inequality to the other.

At the final stage of this work, we notice some recent works by Mijena and Nane \cite{MijenaNane14Int,MijenaNane14ST}, who have also studied this equation
in a more general setting where the Laplacian is replaced by
$-(-\Delta)^{\alpha/2}$ and the space dimension can be any $d<\alpha\: (2\wedge \beta^{-1})$.
When $\beta\in \:]0,1[\:$, $\alpha=2$ and $d=1$, they obtain the same rate as in \eqref{E:WI}.
The main differences of our work from \cite{MijenaNane14Int,MijenaNane14ST}
include:
(1) Our initial data are more general (measures), which entails more calculations;
(2) We cover the case $\beta\in \:]1,2[\:$, to which most efforts in Section \ref{S:Green-Basic} are contributed;
(3) We derive both upper and lower moment bounds, which can be handy for proving many other results;
(4) We prove the weak intermittency of type I for the slow diffusion equations, thanks to our {\em lower} bound on the second moment.

These studies are far from being conclusive. Many aspects can be improved,
such as the H\"older regularity for measure-valued initial data and for fast diffusion equations,
full intermittency for both slow and fast diffusion equations, etc.
Finally, one interesting question is whether the sample-path comparison principle holds for
the slow diffusion equations; see the recent work \cite{ChenKim14Comparison} for the SFHE \eqref{E:FracHeat} and references therein.

\bigskip
This paper is structured as follows.
We first introduce some notation in Section \ref{S:Note}. The main results are stated in Section \ref{S:Main}.
In Section \ref{S:Green-Basic}, we prove some useful properties of the Green functions.
Section \ref{S:K-bounds} gives a general framework on calculating the function $\calK(t,x)$, based on which Theorem \ref{T:K-bounds} is proved.
The proof of the existence and uniqueness results with moment estimates, i.e., Theorem \ref{T:ExUni},
is presented in Section \ref{S:ExUni}.

\paragraph{Acknowledgements}
The author thanks Erkan Nane for pointing out that the classical Duhamel principle fails and one should use
the time-fractional Duhamel principle \cite{UmarovSaydamatov06} as in \cite{MijenaNane14Int,MijenaNane14ST}.
The author thanks Davar Khoshnevisan for some useful comments.

\section{Some preliminaries and notation}
\label{S:Note}
Recall that the Green functions $G_\beta(t,x)$ solve \eqref{E:Homo-Eq}.
Note that in \cite{MainardiEtc01Fundamental},
the fundamental solution is defined with
the initial conditions $u(0,x) = \delta_0(x)$ and
$\frac{\partial}{\partial t} u(0,x) =0$ for all $\beta\in \;]0,2]$.
Let $G_\beta^*(t,x)$, which is also called the Green function, be the solution to \eqref{E:Homo-Eq} subject to the
initial data
\[
u(0,x) = \delta_0(x)\quad\text{and}\quad \frac{\partial }{\partial t} u(0,x) =0.
\]
Here are some special cases. If $\beta=1$, then $G_\beta(t,x)$
reduces to the heat kernel function, i.e.,
\begin{align}\label{E:HeatG}
G_1(t,x)=\frac{1}{\sqrt{4 t}} \exp\left(-\frac{x^2}{4t}\right)\;,\quad\text{for
$(t,x)\in\R_+\times\R$.}
\end{align}
If $\beta=2$, then $G_\beta(t,x)$ and $G_\beta^*(t,x)$
reduce to the heat kernel functions, i.e.,
\begin{align}\label{E:WaveG}
 G_2(t,x) = \frac{1}{2} \Indt{|x|\le t},\quad \text{and}\quad
G_2^*(t,x) = \frac{1}{2}\left(\delta_{t}(x)+\delta_{-t}(x)\right).
\end{align}

For $\mu$ and $\nu\in \calM_T^\beta(\R)$, the solution to the following
homogeneous equation
\begin{align*}
 \begin{cases}
  \left(\displaystyle \DCap -
\frac{\partial^2}{\partial x^2} \right) u(t,x) =
0,& t\in \R_+^*,\:
x\in\R,
\\[1em]
\displaystyle
u(0,\cdot) = \mu(\cdot)\;,&\text{if $0<\beta\le 1$\;,}\\[0.4em]
\displaystyle
u(0,\cdot) = \mu(\cdot)\;,\quad \frac{\partial }{\partial t} u(0,\cdot) =
\nu(\cdot)\;,&\text{if $1<\beta< 2$\;,}
\end{cases}
\end{align*}
will always be denoted by $J_0(t,x)$, which is equal to
\begin{align}\label{E:J0}
J_0(t,x): =
\begin{cases}
\displaystyle
 \int_\R \mu(\ud y)\: G_\beta(t,x-y)\;,& \text{if
$0<\beta\le 1$\;,}\\[0.8em]
\displaystyle
\int_\R \nu(\ud y)\: G_\beta(t,x-y)+\int_\R \mu(\ud y)\:  G_\beta^*(t,x-y)\;,&
\text{if
$1<\beta< 2$\;.}
\end{cases}
\end{align}

\begin{remark}
For the slow diffusion equations ($0<\beta\le 1$), the Green function
$G_\beta(t,x)$ is the same as the function $G_{\alpha,\beta}^\theta(x,t)$ in
\cite[Section 3]{MainardiEtc01Fundamental} with $\alpha=2$ and $\theta=0$.
For the fast diffusion equations ($1<\beta<2$), our function $G_\beta^*(t,x)$
corresponds to the function $G_{\alpha,\beta}^\theta(x,t)$ in
\cite[Section 3]{MainardiEtc01Fundamental}.
In these two cases,  the Green functions $G_\beta(t,x)$ and $G_\beta^*(t,x)$, and their properties are  mostly
known; see \cite{MainardiEtc01Fundamental} and \cite[Appendix F]{Mainardi10Book}.
However, for the fast diffusion equations, the Green function $G_\beta(t,x)$ and
its properties need to be proved, which is
done in Lemma \ref{L:Green} below.
\end{remark}

Let $W=\left\{W_t(A):A\in\calB_b\left(\R\right),t\ge 0 \right\}$
be a space-time white noise
defined on a complete probability space $(\Omega,\calF,P)$, where
$\calB_b\left(\R\right)$ is the collection of Borel sets with finite Lebesgue measure.
Let
\[
\calF_t = \sigma\left(W_s(A):0\le s\le
t,A\in\calB_b\left(\R\right)\right)\vee
\calN,\quad t\ge 0,
\]
be the natural filtration augmented by the $\sigma$-field $\calN$ generated
by all $P$-null sets in $\calF$.
We use $\Norm{\cdot}_p$ to denote the $L^p(\Omega)$-norm ($p\ge 1$).
In this setup, $W$ becomes a worthy martingale measure in the sense of Walsh
\cite{Walsh86}, and $\iint_{[0,t]\times\R}X(s,y) W(\ud s,\ud y)$ is well-defined in this reference for a suitable
class of random fields $\left\{X(s,y),\; (s,y)\in\R_+\times\R\right\}$.
\bigskip

Recall that the rigorous meaning of the spde \eqref{E:tFracHt} is in the integral form \eqref{E:WalshSI}.

\begin{definition}\label{DF:Solution}
A process $u=\left(u(t,x),\:(t,x)\in\R_+^*\times\R\right)$  is called a {\it random field solution} to
\eqref{E:tFracHt} if
\begin{enumerate}[(1)]
 \item $u$ is adapted, i.e., for all $(t,x)\in\R_+^*\times\R$, $u(t,x)$ is $\calF_t$-measurable;
\item $u$ is jointly measurable with respect to
$\calB\left(\R_+^*\times\R\right)\times\calF$;
\item $\left(G_\beta^2 \star \Norm{\rho(u)}_2^2\right)(t,x)<+\infty$
for all $(t,x)\in\R_+^*\times\R$, where $\star$ is the convolution in both
space and time variables. Moreover the function $(t,x)\mapsto I(t,x)$ mapping $\R_+^*\times\R$
into $L^2(\Omega)$ is continuous;
\item $u$ satisfies \eqref{E:WalshSI} a.s.,for all $(t,x)\in\R_+^*\times\R$.
\end{enumerate}
\end{definition}

Assume that the function $\rho:\R\mapsto \R$ is globally Lipschitz
continuous with Lipschitz constant $\LIP_\rho>0$.
We need some growth conditions on $\rho$:
assume that for some constants $\Lip_\rho>0$ and
$\Vip \ge 0$,
\begin{align}\label{E:LinGrow}
|\rho(x)|^2 \le \Lip_\rho^2 \left(\Vip^2 +x^2\right),\qquad \text{for
all $x\in\R$}.
\end{align}
Sometimes we need a lower bound on $\rho(x)$:
assume that for
some constants $\lip_\rho>0$ and $\vip \ge 0$,
\begin{align}\label{E:lingrow}
|\rho(x)|^2\ge \lip_{\rho}^2\left(\vip^2+x^2\right),\qquad \text{for all $x\in\R$}\;.
\end{align}

For all $(t,x)\in\R_+^*\times \R$, $n\in\bbN$ and $\lambda\in\R$, define
\begin{align}
\notag
\calL_0\left(t,x;\lambda\right) &:= \lambda^2 G_\beta^2(t,x)  \\
\label{E:Ln}
\calL_n\left(t,x;\lambda\right)&:=
\left(\calL_0\star \cdots\star\calL_0\right)(t,x),\quad\text{for $n\ge 1$, ($n$ convolutions),}\\
\calK\left(t,x;\lambda\right)&:= \sum_{n=0}^\infty
\calL_n\left(t,x;\lambda\right).
\label{E:K}
\end{align}
We will use the following conventions to the kernel functions $\calK(t,x;\lambda)$:
\begin{align*}
 \calK(t,x) &:= \calK(t,x;\lambda),&
\overline{\calK}(t,x) &:= \calK\left(t,x;\Lip_\rho\right),\\
\underline{\calK}(t,x) &:= \calK\left(t,x;\lip_\rho\right),&
\widehat{\calK}_p(t,x) &:= \calK\left(t,x; 4\sqrt{p} \Lip_\rho\right),\quad\text{for $p\ge 2$}\:.
\end{align*}

\section{Main results}
\label{S:Main}
Our first theorem is about the existence, uniqueness and moment estimates of the solutions to \eqref{E:tFracHt}.
It possesses a general form as \cite[Theorem 2.4]{ChenDalang13Heat}, \cite[Theorem 2.3]{ChenDalang14Wave}, and
\cite[Theorem 3.1]{ChenDalang14FracHeat}.

\begin{theorem}[Existence,uniqueness and moments]\label{T:ExUni}
Suppose that\\
(i) $0<\beta< 2$;\\
(ii) The function $\rho$ is Lipschitz continuous and satisfies the growth
condition \eqref{E:LinGrow};\\
(iii) The initial data are such that $\mu\in\calM_T^\beta\left(\R\right)$ if $\beta \in\:]0,1]$,
and
$\mu, \: \nu\in\calM_T^\beta\left(\R\right)$ if $\beta\in \:]1,2[\:$.\\
Then the SPDE \eqref{E:tFracHt} has a unique (in the sense of versions) random field solution $\{u(t,x): (t,x)\in\R_+^* \times
\R \}$. Moreover, the following statements are true:\\
(1) $(t,x)\mapsto u(t,x)$ is $L^p(\Omega)$-continuous for all integers $p\ge
2$; \\
(2) For all even integers $p\ge 2$, all $t>0$ and
$x,y\in\R$,
\begin{align}\label{E:MomUp}
 \Norm{u(t,x)}_p^2 \le
\begin{cases}
 J_0^2(t,x) + \left(\left[\Vip^2+J_0^2\right] \star \overline{\calK} \right)
(t,x),& \text{if $p=2$}\;,\cr
2J_0^2(t,x) + \left(\left[\Vip^2+2J_0^2\right] \star \widehat{\calK}_p \right)
(t,x),& \text{if $p>2$}\;;
\end{cases}
\end{align}
(3) If $\rho$ satisfies \eqref{E:lingrow}, then
for  all $t>0$ and $x,y\in\R$,
\begin{align}
\label{E:SecMom-Lower}
\Norm{u(t,x)}_2^2 \ge J_0^2(t,x) + \left(\left(\vip^2+J_0^2\right) \star
\underline{\calK} \right)
(t,x)\;.
\end{align}
\end{theorem}

The following Theorem \ref{T:Holder} gives the H\"older continuity of the solution for the slow diffusion equations.
We cannot prove the H\"older regularity for the fast diffusion equations due to the less precise results in Proposition \ref{P:G-FD} than those in Proposition
\ref{P:G-SD}.

\begin{theorem}
\label{T:Holder}
Suppose that $\beta\in \:]0, 1]$. If $\mu(\ud x) = f(x)\ud x$
with $f \in
L^\infty\left(\R\right)$, then
\begin{align}\label{E:LpBdd}
\sup_{(t,x)\in [0,T]\times\R}
\Norm{u(t,x)}_p^2 <+\infty, \quad\text{for all $T\ge 0$ and $p\ge 2$.}
\end{align}
Moreover, we have
\begin{align}\label{E:Holder-I}
I(\cdot,\circ)\in
C_{\frac{2-\beta}{4}-,\frac{1}{2}-}\left(\R_+\times\R\right)\;,\quad\text{a.s.,}
\end{align}
and \eqref{E:Holder-u} holds.
If $f$ is bounded and $\alpha$-H\"older continuous
($\alpha\in \;]0,1[\:$), then \eqref{E:HHolder-u} holds.
\end{theorem}
\begin{proof}
The bound  \eqref{E:LpBdd} is a simple consequence of \eqref{E:MomUp}.
The proof of \eqref{E:Holder-I} is straightforward under \eqref{E:LpBdd} (see
\cite[Remark \myRef{4.6}{RH:BddHolder}]{ChenDalang13Holder}).
The rest parts are due to Lemma \ref{P:J0LocalLip}.
\end{proof}

In only very few cases, one can derive explicit form for $\calK(t,x)$.
A first case is when $\beta=1$; see Example \ref{Ex:SHE}.
A second case is given in Example \ref{Ex:DoubleSHE}.
A third case is when $\beta =2$:
\[
\calK^{\mbox{\scriptsize wave}}(t,x;\lambda) = \frac{\lambda^2}{4}
I_0\left(\sqrt{\frac{\lambda^2((\kappa t)^2-x^2)}{2\kappa}}\right) \Indt{|x|\le
\kappa t},
\]
where $I_0(x)$ is the {\it modified Bessel function of the first kind}
of order $0$; see \cite{ChenDalang14Wave}.
Hence, in order to use the moment bounds in \eqref{E:MomUp} and \eqref{E:SecMom-Lower},
we need some good estimates on the kernel function $\calK(t,x)$.
For this purpose, we define some {\it reference kernel functions}:
\begin{align}\label{E:calG}
\calG_\beta(t,x) :=
\begin{cases}\displaystyle
\frac{1}{2\:
t^{\beta/2}}\exp\left(-\frac{|x|}{t^{\beta/2}}\right)
&\text{if $0<\beta<1$\;.}\\[0.6em]\displaystyle
\frac{1}{\sqrt{4\pi t^\beta
}}\exp\left(-\frac{x^2}{4t^\beta}\right)& \text{if $1\le \beta<2$\;.}
\end{cases}
\end{align}
Note that when $1\le \beta <2$, $\calG_\beta(t,x)=G_1\left(t^\beta,x\right)$.
For convenience, when $0<\beta<1$, denote
\begin{align}\label{E:Geb}
\calG_{e,\beta}(t,x) :=\calG_{\beta}(t,x),
\end{align}
where the subscription ``$e$'' refers to the exponential function§.
Clearly, $\calG_\beta(t,x)$ is nonnegative and
$\int_\R \ud x\: \calG_\beta(t,x) =1$. For $0<\beta <1$, define
\begin{align}\label{E:calG-L}
\ubar{\calG}_\beta(t,x):= G_1\left(t^\beta,x\right
) = \frac{1}{\sqrt{4\pi
t^\beta}} \exp\left(-\frac{x^2}{4\:t^\beta}\right).
\end{align}
We need some constants:
\begin{align}\label{E:HatC}
\widehat{C}_\beta :=\frac{2^{\beta/2 }}{2^{\beta/2 }-1} \exp\left(-\frac{1}{2^{\beta/2}}\right),\quad\text{for $\beta>0$,}
\end{align}
and
\begin{align}\label{E:TildeC}
\widetilde{C}_\beta :=
\begin{cases}
\widehat{C}_\beta
& \text{if $0<\beta<1$\;,}\\[0.5em]
 2^{\frac{\beta-1}{2}}
& \text{if $1\le \beta<2$\;.}
\end{cases}
\end{align}
\begin{remark}
The constant $\widehat{C}_\beta$ as a function of $\beta\in \:]0,2]$ is decreasing with $\widehat{C}_2=2 e^{-1/2}\approx 1.21306$, $\widehat{C}_1=\left(2+\sqrt{2}\right)
e^{-\frac{1}{\sqrt{2}}}\approx 1.68344$, and $\lim_{\beta\rightarrow 0_+} \widehat{C}_\beta = \infty$.
\end{remark}
Define
\begin{align}\label{E_:G-Bdd}
\Psi_\beta:=\sup_{x\in\R}
\frac{G_\beta^2(1,x)}{\calG_\beta(1,x)}<+\infty\;,\quad\text{for $0<\beta<2$,}
\end{align}
and
\begin{align}\label{E_:G-Bdd-L}
\ubar{\Psi}_\beta :=\inf_{x\in\R}
\frac{G_\beta^2(1,x)}{\ubar{\calG}_\beta(1,x)}>0,
\quad\text{for $0<\beta<1$.}
\end{align}
Proposition \ref{P:ST-Con} below shows that
$\Psi_\beta<+\infty$ and $\ubar{\Psi}_\beta>0$.

\begin{theorem}
\label{T:K-bounds}
Fix $\lambda>0$.
(1) For $\beta\in\:]0,2[\:$, there is a finite constant $C:=C(\beta,\lambda)$ such that
\begin{align}\label{E:calK-U}
\calK(t,x;\lambda)\le \frac{C}{t^\sigma} \; \calG_\beta(t,x) \left(1+ t^\sigma
\exp\left(\Upsilon t\right)\right),
\end{align}
where
\begin{align}\label{E:SigUps}
\sigma=\beta/2 +2(1-\Ceil{\beta})\quad
\text{and}\quad
\Upsilon= \left(\lambda^2\: \Psi_\beta \;\widetilde{C}_\beta\:
\Gamma(1-\sigma)\right)^{\frac{1}{1-\sigma}}.
\end{align}
(2) For $\beta\in\:]0,1[\:$, there is a constant $\ubar{C}
:=\ubar{C}(\beta,\lambda)>0$
such that
\begin{align}
\label{E:calK-L}
\calK(t,x;\lambda)\ge \ubar{C} \; \ubar{\calG}_\beta(t,x)
 \exp\left(\ubar{\Upsilon}\: t\right),
\end{align}
where
\[
\ubar{\sigma}=\beta/2-1\quad\text{and}\quad
\ubar{\Upsilon}=\left(2^{-1/2} \: \lambda^2\: \ubar{\Psi}_\beta \:
\Gamma(1-\ubar{\sigma})\right)^{\frac{1}{1-\ubar{\sigma}}}.
\]
\end{theorem}
\begin{proof}
Apply Proposition \ref{P:ST-Con} below with $\lambda G_\beta(t,x)$. Note that introducing the factor $\lambda$
changes the constants $C_0$ and $\ubar{C}_0$ by a factor $\lambda^2$.
\end{proof}

The last set of results are the weak intermittency and the bounds in \eqref{E:WI}.

\begin{theorem}[Weak intermittency of type I for slow diffusion equations]
\label{T:Weak-S}
Suppose that $\beta \in\;]0,1[\:$ and $\mu(\ud x) = c \:\ud x$. If $\rho$
satisfies \eqref{E:LinGrow} and $|c|+|\Vip|\ne 0$,
then
\[
\overline{m}_p\le \frac{1}{2}\left[2^4 \Lip_\rho^2\: \widehat{C}_\beta  \Psi_\beta \:
\Gamma(1-\beta/2)\right]^{\frac{2}{2-\beta}}\:
p^{\frac{4-\beta}{2-\beta}},\quad\text{for all $p\ge 2$ even.}
\]
If $\rho$ satisfies  \eqref{E:lingrow} and $|c|+|\vip|\ne 0$, then the solution is
{\em weakly intermittent of type I}:
\[
\underline{m}_p \ge\frac{p}{2}
\left(2^{-1/2}\lip_\rho^2 \: \ubar{\Psi}_\beta \:
\Gamma(2-\beta/2)\right)^{\frac{2}{4-\beta}}
 ,\quad\text{for all $p\ge 2$.}
\]
\end{theorem}
\begin{proof}
Clearly, in this case,
$J_0(t,x) = c$. Hence, by \eqref{E:MomUp} and \eqref{E:calK-U},
\[
\Norm{u(t,x)}_p^2 \le c^2 + \frac{C}{t^\sigma}\left(\Vip^2+2 c^2\right)
\left(1+t^\sigma \exp\left(\left[2^4 \Lip_\rho^2 \:\widehat{C}_\beta \Psi_\beta \:
\Gamma(1-\sigma)\right]^{\frac{1}{1-\sigma}} p^{\frac{1}{1-\sigma}}\:
t\right)\right),
\]
with $\sigma=\beta/2$. Then increase the power by a factor $p/2$. This proves the upper
bounds. As for the lower bound, by \eqref{E:SecMom-Lower} and \eqref{E:calK-L},
\[
\Norm{u(t,x)}_p^2\ge
\Norm{u(t,x)}_2^2
\ge
c^2 + \ubar{C}  \left(\vip^2 + c^2\right)
 \exp\left( \left[2^{-1/2}\lip_\rho^2
\ubar{\Psi}_\beta \: \Gamma(1-\ubar{\sigma})\right
]^\frac{1}{1-\ubar{\sigma}} t
\right)
\]
with $\ubar{\sigma}= \frac{\beta}{2}-1$. This
completes the proof.
\end{proof}

\begin{theorem}[Weak intermittency of type II for fast diffusion equations]
\label{T:Weak-F}
Suppose that $\beta\in \:]1,2[\:$, $\mu(\ud x) = c\: \ud x$ and $\nu(\ud x) = c'\: \ud x$.
If $\rho$ satisfies \eqref{E:LinGrow} and $|c|+|c'|+|\Vip|\ne 0$, then
\[
\overline{m}_p\le
\frac{1}{2}\left[2^{9/2} \Lip_\rho^2\:
\Psi_\beta \:
\Gamma(3-\beta/2)\right]^{\frac{2}{6-\beta}}\:
p^{\frac{8-\beta}{6-\beta}},\quad\text{for all $p\ge 2$
even.}
\]
\end{theorem}
\begin{proof}
By Lemma \ref{L:Green} (iii), $J_0(t,x) = c \:t+ c'$.
The condition $|c|+|c'|\ne 0$ implies $J_0(t,x)\ne 0$ for large $t$.
Hence, by \eqref{E:MomUp} and \eqref{E:calK-U},
\[
\Norm{u(t,x)}_p^2 \le (c\:t+c')^2 +  \frac{C}{t^\sigma}\left(\Vip^2+2
(c\:t+c')^2\right)
\left(1+t^\sigma \exp\left(\left[2^{4}\widetilde{C}_\beta \Lip_\rho^2 \Psi_\beta \:
\Gamma(1-\sigma)\right]^{\frac{1}{1-\sigma}}
p^{\frac{1}{1-\sigma}}\: t\right)\right)
\]
with $\sigma=\frac{\beta}{2}-2$. Then increase the power by a factor $p/2$ and use the fact that $\widetilde{C}_\beta\le \sqrt{2}$.
\end{proof}


\section{Some properties of the Green functions}
\label{S:Green-Basic}
We need some special functions. The following two-parameter
{\it Mittag-Leffler function}
\begin{align}\label{E:Mittag-Leffler}
E_{\alpha,\beta}(z) := \sum_{k=0}^{\infty}
\frac{z^k}{\Gamma(\alpha k+\beta)},
\qquad \alpha>0,\;\beta> 0\;,
\end{align}
is a generalization of exponential function, $E_{1,1}(z)=e^z$; see, e.g.,
\cite[Section 1.2]{Podlubny99FDE}. Another special case\footnote{Proof of \eqref{E:ML-HH}.
By \cite[41:6:6]{oldham2008atlas},
$\frac{1}{\sqrt{\pi x^2 }}+
e^{x^2}\Erfc(-x) = \frac{1}{x^2}\sum_{n=1}^\infty
\frac{x^n}{\Gamma(n/2)}$. Then apply \eqref{E:Efrom1} below.\myEnd} is
\begin{align}\label{E:ML-HH}
 E_{1/2,1/2}(x) =  \frac{1}{\sqrt{\pi}}+
x\:e^{x^2} \Erfc(-x),\quad\text{for $x\ge 0$},
\end{align}
where $\Erf(x)=\frac{2}{\sqrt{x}}\int_0^x\ud y \:e^{-y^2}$ is
the {\it error function} and $\Erfc(x)=1-\Erf(x)$ is the
{\it complementary  error function}.
We will use the convention that $E_{\alpha}(z)=E_{\alpha,1}(z)$.
A function is called {\it completely monotonic}
if  $(-1)^n f^{(n)}(x)\ge 0$ for $n =0,1,2,\dots$;
 see \cite[Definition 4, on
p. 108]{Widder41LaplaceTr}.
An important fact \cite{Schneider96CM} that we
are going  to use is that
\begin{align}\label{E:E-CM}
 \text{$x\in\R_+\mapsto E_{\alpha,\beta}(-x)$ is
completely monotonic}
 \quad\Longleftrightarrow\quad
0<\alpha\le 1\wedge \beta.
\end{align}

Let $W_{\lambda,\mu}(z)$ be the {\it two-parameter
 Wright function} of order $\lambda$ defined as
follows:
\[
W_{\lambda,\mu}(z) := \sum_{n=0}^\infty \frac{z^n}
{n!\; \Gamma\left(\lambda
n+ \mu\right)}\;,\qquad\text{for $\lambda>-1$, $\mu\in\bbC$ and $z\in\bbC$\:;}
\]
see, e.g., \cite[Appendix F]{Mainardi10Book} and
references therein.
We define the {\it two-parameter Mainardi functions} of order
$\lambda\in [0,1[\;$ by
\begin{align}\label{E:2p-Mainardi}
M_{\lambda,\mu} (z) :=W_{-\lambda,\mu-\lambda}(-z)
=\sum_{n=0}^\infty \frac{(-z)^n}{n! \;\Gamma\left(
\mu-(n+1)\lambda\right)}
\;,\quad\text{for $\mu\in\bbC$ and $z\in\bbC$\:,}
\end{align}
and we will use the convention that
$M_{\lambda}(z)=M_{\lambda,1}(z)$.
In particular,
$M_{1/2}(z)=\frac{1}{\sqrt{\pi}}\exp\left(-z^2/4\right)$.
The one-parameter Mainardi functions $M_\lambda(z)$ are used by
Mainardi, {\it et al} in \cite{MainardiEtc01Fundamental,Mainardi10Book}.
This two-parameter extension is necessary for the Green function $G_\beta(t,x)$ of the fast diffusions.

\begin{figure}[h!tbp]
\centering
\subfloat[$G_\beta(1,x)$]{
\includegraphics[scale=0.7]{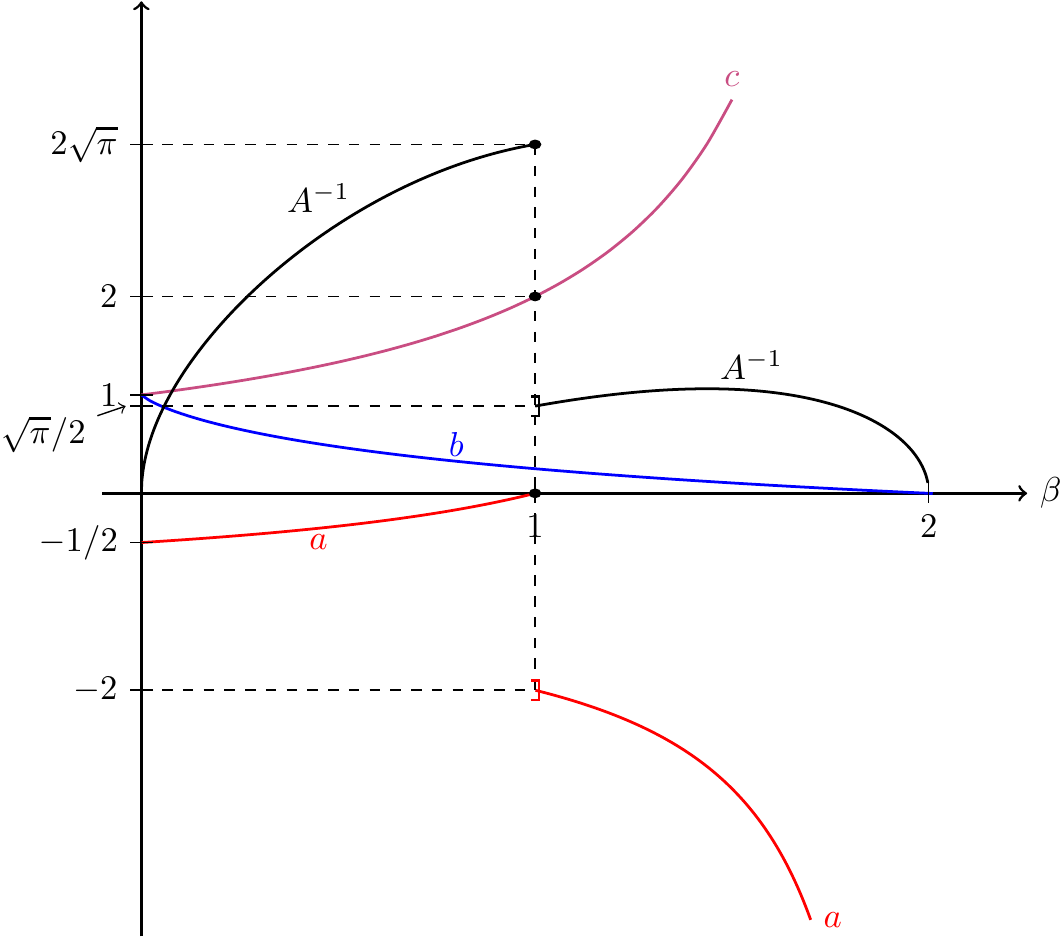}}
\subfloat[$G_\beta^*(1,x)$]{
\raisebox{24.5mm}{\includegraphics[scale=0.7]{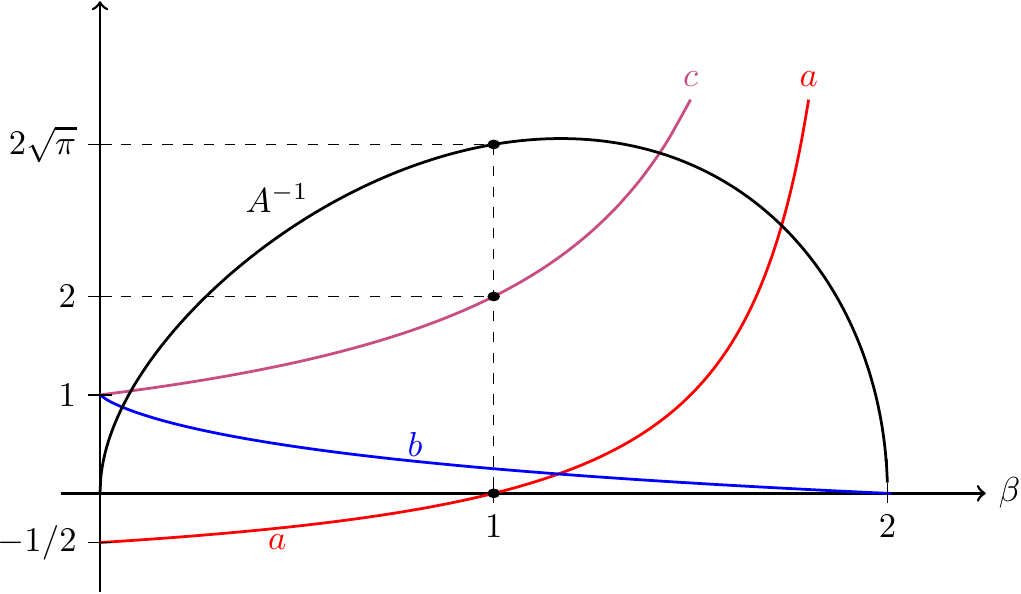}}}
\caption{The parameters of the asymptotics of the
functions $G_\beta(1,x)$ and $G_\beta^*(1,x)$.}
\label{F6:Asy-Parameters}
\end{figure}

\begin{lemma}[Properties of the Green functions $G_\beta(t,x)$ and $G_\beta^*(t,x)$]
\label{L:Green}
For $\beta\in \;]0,2[\;$, the following properties
 hold:\\
(i) The Green function $G_\beta(t,x)$ has the following explicit form
\begin{align}  \label{E:Green}
\hspace{-1.2em}
 G_{\beta}(t,x) &=
\frac{t^{\Ceil{\beta}-1-\beta/2}}{2} \;
M_{\beta/2,\Ceil{\beta}}\left(\frac{|x|}{t^{\beta/2}}\right)=
\begin{cases}
\displaystyle
\frac{t^{-\beta/2}}{2}
M_{\beta/2}\left(\frac{|x|}{t^{\beta/2}}\right)\;,
& \text{if
$0<\beta\le 1$},\\[1em]
\displaystyle
\frac{t^{1-\beta/2}}{2}
M_{\beta/2,2}\left(\frac{|x|}{t^{\beta/2}}\right)\;,& \text{if
$1<\beta< 2$}.
\end{cases}
\end{align}
The function $G_\beta^*(t,x)$ has the same form
as \eqref{E:Green} except that
all $\Ceil{\beta}$'s in \eqref{E:Green} should be
replaced by $1$, i.e.,
\[
G_\beta^*(t,x) = \frac{t^{-\beta/2}}{2}
M_{\beta/2}\left(\frac{|x|}{t^{\beta/2}}\right), \quad \text{for $1<\beta<2$.}
\]
(ii) $G_\beta(t,x)$ has the following scaling property:
\begin{align}\label{E:Scaling}
 G_\beta(t,x) =
t^{\Ceil{\beta}-1-\beta/2}
G_\beta\left(1,\frac{x}{t^{\beta/2}}\right).
\end{align}
The scaling property of $G_\beta^*(t,x)$ is the
same as \eqref{E:Scaling}
except that the $\Ceil{\beta}$ in \eqref{E:Scaling} should be replaced by $1$.\\
(iii) For any $t>0$ fixed, both functions
$x\mapsto G_\beta(t,x)$ and $x\mapsto G_\beta^*(t,
x)$ are symmetric and
nonnegative, i.e., $G_\beta(t,x)=G_\beta(t,-x)\ge 0$ and
$G_\beta^*(t,x)=G_\beta^*(t,-x)\ge 0$, for all $x\in\R$. Moreover
\begin{align}\label{E:IntG}
\int_\R \ud x\: G_\beta(t,x) =
t^{\Ceil{\beta}-1}\quad\text{and}\quad\int_\R\ud x\: G_\beta^*(t,x) =1\;.
\end{align}
In particular, the functions
$x\mapsto G_\beta(t,x)$ with $\beta\in\:]0,1]$ and $x\mapsto
G_\beta^*(t,x)$ are probability densities.\\
(iv) $G_\beta(1,x)$ has the following asymptotic property:
\begin{align}\label{E:AsymG}
 G_\beta(1,x) =\frac{1}{2}M_{\beta/2,\Ceil{\beta}}
\left(|x|\right) \approx
A \:|x|^a\:  e^{-b |x|^c},\quad\text{as $|x|\rightarrow +\infty$}\;,
\end{align}
where
\begin{gather}\label{E:Asym-A}
A=\left(2\pi
(2-\beta)2^{\frac{\beta+4\left(1-\Ceil{\beta}\right)}{2-\beta}}
\;\beta^{\frac{2\left(1+\beta-2\Ceil{\beta}\right)
}{\beta-2}}\right)^{-1/2}
\ge\frac{1}{\sqrt{4\pi}}\;,\\
\label{E:Asym-a}
a= \frac{1+\beta-2\Ceil{\beta}}{2-\beta}\le 0\;,\\
\label{E:Asym-c}
c=\frac{2}{2-\beta}\;>1\;,\\
\label{E:Asym-b}
b=\left(2-\beta\right)2^{-2/\left(2-\beta\right)}
\beta^{\beta/\left(2-\beta\right)}\in \;]0,1[\;.
\end{gather}
$G_\beta^*(1,x)$ has the same asymptotic property except that all
$\Ceil{\beta}$'s in \eqref{E:AsymG}, \eqref{E:Asym-A}
and \eqref{E:Asym-a} should be replaced by $1$, the range of $a$ is
$\;]-1/2,+\infty[\;$ and the range of $A$ is $[1/\sqrt{4\pi},+\infty[\:$.
See Figure \ref{F6:Asy-Parameters} for the plots of these parameters as
functions of $\beta$.\\
(v) $G_\beta(t,x)$ satisfies the following moment formula:
\begin{align}\label{E:G-Mom}
 \int_\R \ud x\: |x|^{a} G_\beta(t,x) =
\frac{\Gamma\left(a+1\right)}{\Gamma\left(\frac{a\:
\beta}{2}+\Ceil{\beta}\right)} \: t^{a \beta/2+\Ceil{\beta}-1}\;,\quad
\text{for $a>-1$ and $t\ge 0$}\:.
\end{align}
The moment formula for $G_\beta^*(t,x)$ is the
same as \eqref{E:G-Mom}
except that all $\Ceil{\beta}$'s should be replace
d by $1$.\\
(vi) The Fourier transform of the Green function
$G_\beta(t,x)$ is
\begin{align}\label{E:G-Fourier}
\int_\R  \ud x\: e^{-i \xi x} G_\beta(t,x)= t^{\Ceil{\beta}-1}\;
E_{\beta,\Ceil{\beta}}\left(-t^\beta
\xi^2\right)\;,\quad\text{for $t>0$ and $\xi\in\R$\:.}
\end{align}
The Fourier transform of $G_\beta^*(t,x)$ is the
same as \eqref{E:G-Fourier}
except that all $\Ceil{\beta}$'s in \eqref{E:G-Fourier} should be replaced by
$1$.\\
(vii) the Laplace transform of the function $\R_+\ni x\mapsto G_\beta(1,x)$ is
\begin{align}\label{E:G-Laplace}
 \int_0^\infty\ud x\:  e^{- z x} G_\beta(1,x) =\frac{1}{2}
E_{\beta/2,\Ceil{\beta}}(-z)\;,\quad\text{for all $z\in\bbC$} \:.
\end{align}
The Laplace transform of $\;\R_+\ni x\mapsto G_\beta^*(1,x)$ is the same as
\eqref{E:G-Laplace} except that the $\Ceil{\beta}$ in \eqref{E:G-Laplace} should be replaced by $1$.\\
(viii) The function $x\mapsto G_\beta(t,x)$ attains
 its maximum value
at $x=0$:
\begin{align}\label{E:G-Max}
\sup_{x\in\R} G_\beta(t,x) = G_\beta(t,0) =
\frac{t^{\Ceil{\beta}-1-\beta/2}}{2}
\Gamma\left(\Ceil{\beta}-\frac{\beta}{2}
\right)^ { -1} \;;
\end{align}
The function  $x\mapsto G_\beta^*(t,x)$ attains two symmetric maximums that move apart from the origin with time.\\
(ix) The function $x\mapsto G_\beta(t,x)$ is continuous at $x=0$ but in general not differentiable there.
Its $n$-th derivatives are equal to
\begin{align}\label{E:G-dx}
\frac{\partial^n}{\partial x^n} G_\beta(t,x)
=
\begin{cases}
 \displaystyle
\frac{(-1)^n t^{\Ceil{\beta}-1-(n+1)\beta/2}}{2}
M_{\beta/2,\Ceil{\beta}-n\beta/2}\left(\frac{x}{
t^{\beta/2}}\right)&
\text{if $x>0$\;,}\\[0.8em]
\displaystyle
\frac{t^{\Ceil{\beta}-1-(n+1)\beta/2}}{2}
M_{\beta/2,\Ceil{\beta}-n\beta/2}\left(-\frac{x}{
t^{\beta/2}}\right)&
\text{if $x<0$\;.}
\end{cases}
\end{align}
\end{lemma}
\begin{proof}
Denote
\[
G_\beta^\star(t,x):=\begin{cases}
G_\beta(t,x)& \text{if $\beta\in \:]0,1]$},\cr
G_\beta^*(t,x)& \text{if $\beta\in \:]1,2[$}\:.
\end{cases}
\]
All these properties for $G_\beta^\star(t,x)$ can be found
in \cite{MainardiEtc01Fundamental} and \cite[Appendix F]{Mainardi10Book}.
The expression \eqref{E:Green} for $G_\beta^\star(t,x)$ can be found in
\cite[(4.23)]{MainardiEtc01Fundamental}.
The scaling property \eqref{E:Scaling} for $G_\beta^\star(t,x)$ can be found in
\cite[(3.7)]{MainardiEtc01Fundamental}.
The asymptotic property of $G_\beta^\star(t,x)$ can
be found in \cite[(4.29), (4.30)]{MainardiEtc01Fundamental}.
The moment formula \eqref{E:G-Mom} for $G_\beta^\star(t,x)$ can be found in
\cite[(4.31)]{MainardiEtc01Fundamental}, where one can extend integer $n$ to all $a>-1$.
The Fourier transform of $G_\beta^\star(t,x)$ can be found in
\cite[(4.21)]{MainardiEtc01Fundamental}.
The Laplace transform \eqref{E:G-Laplace} of  $G_\beta^\star(t,x)$ is due to the
Laplace transform of the Wright function of the second kind (see e.g., \cite[(F.25), on p. 248]{Mainardi10Book}):
$W_{-\lambda,\mu}(-x) \stackrel{\calL}{\longrightarrow} E_{\lambda,\lambda+\mu}(-z)$ for $0<\lambda<1$,
which implies
\begin{align}\label{E:M-Lap}
M_{\lambda, \mu}(x) \stackrel{\calL}{\longrightarrow}
E_{\lambda,\mu}(-z)\;,\quad\text{for $0<\lambda<1$ \;.}
\end{align}
The statements in both (iii) and (viii) for  $G_\beta^\star(t,x)$ can be found in
\cite[p. 22]{MainardiEtc01Fundamental}.

It remains to prove properties of the Green functions $G_\beta(t,x)$ with $\beta\in \:]1,2[\:$.
Since the arguments for $G_\beta(t,x)$ with $\beta\in \:]0,1]$ are similar to those for $G_\beta(t,x)$ with $\beta\in \:]1,2[\:$, in the following, we will prove both cases altogether.
We will mostly follow the arguments by Mainardi, {\it et al} in \cite{MainardiEtc01Fundamental}.
Let $\widehat{f}$ and $\widetilde{g}$ denote the Fourier transform in the
space variable and the Laplace transform in the
time variable, respectively.
Apply the Fourier transform on the initial data of \eqref{E:Homo-Eq}:
\[
\begin{cases}
\widehat{G_\beta}(0_+,\xi) = 1\;, &\text{if $0<\beta\le 1$}\;,\\[0.5em]
\displaystyle
\widehat{G_\beta}(0_+,\xi) = 0\:,\quad \frac{\partial }{\partial t}
\widehat{G_\beta}(0_+,\xi) = 1\;,
&\text{if $1<\beta< 2$\:.}
\end{cases}
\]
Apply both the Fourier and the Laplace transforms on the both sides of the main
equation in \eqref{E:Homo-Eq}:
\[
\widetilde{\widehat{G_\beta}}(s,\xi) s^\beta - s^{\beta-\Ceil{\beta}}+\xi^2
\widetilde{\widehat{G_\beta}}(s,\xi)=0\;,
\]
where we have used the equivalent definition of the Caputo fractional
differential operator of order $\beta$ through the Laplace transform (see
\cite[(2.12)]{MainardiEtc01Fundamental}):
\[
\calL\left[\DCap f(t)\right](s) = s^\beta \;\widetilde{f} -
\sum_{k=0}^{m-1} s^{\beta-1-k}\; f^{(k)}(0_+)\;,
\quad \text{if $m-1<\beta\le m$}\;.
\]
Hence,
\[
\widetilde{\widehat{G_\beta}}(s,\xi) = \frac{s^{\beta -\Ceil{\beta}}}{s^\beta
+\xi^2}\;,\quad \text{for $0<\beta< 2$\;.}
\]
By the scaling rules for the Fourier and Laplace transforms, we have that
\begin{align*}
G_\beta(bt,ax) &\stackrel{\calF}{\longrightarrow}
\frac{1}{a} \widehat{G_\beta(b t,\cdot)}(\xi/a)
\stackrel{\calL}{\longrightarrow}
\frac{1}{a b} \widetilde{\widehat{G_\beta}}(s/b,\xi/a)
=
\frac{1}{ab}\frac{\left(\frac{s}{b}\right)^{\beta-
\Ceil{\beta}}}{\left(\frac{s}{b}\right)^{\beta} + \left(\frac{\xi}{a}\right)^2}\\
&=
\frac{1}{ab^{1-\Ceil{\beta}}}\frac{s^{\beta-\Ceil{\beta}}
} { s^{\beta} + \left(\frac{b^{\beta/2}\xi}{a}\right)^2}
\stackrel{\calL^{-1}}{\longrightarrow}
\frac{1}{ab^{1-\Ceil{\beta}}}
\widehat{G_\beta(t,\cdot)}\left(\frac{b^{\beta/2}}
{a} \xi\right)\\
&\stackrel{\calF^{-1}}{\longrightarrow}
b^{-\beta/2+ \Ceil{\beta}-1} G_\beta\left(t,\frac{a}{b^{\beta/2}} x\right)\;,
\end{align*}
which proves the scaling property \eqref{E:Scaling}.
Now use the following Laplace transform (see \cite[(1.80, on p. 21)]{Podlubny99FDE})
\[
\int_0^\infty\ud t\:  e^{-s t}
t^{\alpha k + \beta -1} E_{\alpha,\beta}^{(k)}\left(\pm \lambda
t^{\alpha}\right) = \frac{k!\; s^{\alpha-\beta}}{\left(s^\alpha \mp
\lambda\right)^{k+1}}\;,\qquad \Re(s)>|\lambda|^{1/\alpha}\;,
\]
where $E_{\alpha,\beta}^{(k)}(y) = \frac{\ud^k}{\ud y^k}E_{\alpha,\beta}(y)$.
We see that $\widehat{G_\beta}(t,\xi) = t^{\Ceil{\beta}-1} \;E_{\beta,\Ceil{\beta}}\left(-\xi^2 t^\beta\right)$,
which proves \eqref{E:G-Fourier}.
Then an application of the inverse Fourier transform using
Lemma \ref{L:Mf-Fourier} gives the Green function \eqref{E:Green}.
As a consequence, the function $x\mapsto G_\beta(t,x)$ is symmetric and
\[
\int_\R\ud x\:  G_\beta(t,x) = \widehat{G_\beta}(t,0) = t^{\Ceil{\beta}-1}\;,
\]
which proves \eqref{E:IntG}.
By the scaling property and the symmetry of $x\mapsto G_\beta(t,x)$,
\[
\int_\R \ud x\: |x|^n G_\beta(t,x)
= 2 \int_0^\infty\ud x\: x^n G_\beta(t,x)
=2 \: t^{\frac{n\beta}{2}+\Ceil{\beta}-1}
\int_0^\infty\ud y\: y^n\: G_\beta(1,y).
\]
Then the moment formula \eqref{E:G-Mom} is proved by applying Lemma \ref{L:Mf-Mom}.

The asymptotic property of $G_\beta(1,x)$ is a direct consequence of the
asymptotics of the Wright function (see \cite{Wright34AsymBessel} and also
\cite[(F.3), on p. 238]{Mainardi10Book}): For $0<\lambda<1$, and $\mu\in\R$,
\begin{align}\label{E:Asy-M}
 M_{\lambda,\lambda+\mu}(x)=W_{-\lambda,\mu}(-x) \approx A_0\; Y_\lambda(x)^{1/2-\mu}
\exp\left(-Y_\lambda(x)\right)\;,\quad
\text{as $x\rightarrow+\infty$\;,}
\end{align}
where
\[
Y_\lambda(x)=\left(1-\lambda\right)\lambda^{\frac{\lambda}{1-\lambda}}
x^{\frac{1}{1-\lambda}}\quad\text{and}\quad
A_0=\left(2\pi \left(1-\lambda\right)\right)^{-1/2}\;.
\]

The Laplace transform in \eqref{E:G-Laplace} is proved by
\eqref{E:M-Lap}. Bernstein's theorem on monotone functions (see Theorem
\ref{T:Bernstein}) and \eqref{E:E-CM} prove the positivity of $G_\beta(1,x)$ for $x\ge 0$.
Then by symmetry of $G_\beta(1,x)$,  $G_\beta(1,x)\ge 0$ for all $x\in\R$.
By \eqref{E:G-Laplace} and the property of the Laplace transform, we see that
\begin{align*}
\calL\left[\frac{\ud }{\ud x} G_\beta(1,\cdot) \right](z) &=
\frac{z}{2} E_{\beta/2,\Ceil{\beta}}(-z) -
G_\beta(1,0)=-E_{\beta/2,\Ceil{\beta} -\beta/2 }(-z)\;,
\end{align*}
where we have also used the fact that
$G_\beta(1,0) =\left[2\Gamma\left(\Ceil{\beta}-\beta/2\right)\right]^{-1}$
and the recurrence relation of the Mittag-Leffler function
$E_{\alpha,\beta}(z) = \Gamma\left(\beta\right)^{-1} +z E_{\alpha,\alpha+\beta}(z)$.
Notice the function $E_{\beta/2,\Ceil{\beta} -\beta/2 }(-x)$ is complete monotone for
$x\in\R_+$ because $\Ceil{\beta}-\beta/2\ge \beta/2$ and $\beta/2\in\;]0,1]$.
Hence, by the same reason for the positivity of $G_\beta(1,x)$, we can conclude
the non-positivity of $\frac{\ud}{\ud x} G_\beta(1,x)$,
which proves that the global maximum of $G_\beta(t,x)$ is achieved at $x=0$.

As for \eqref{E:G-dx}, by differentiating term-by-term (see also \cite[(F. 8), on p. 239]{Mainardi10Book}),
we see that $\frac{\ud }{\ud z} W_{\lambda,\mu}(z)= W_{\lambda,\lambda+\mu}(z)$,
from which one can easily derive that
\begin{align}\label{E:Der-M}
\frac{\ud^n }{\ud z^n} M_{\lambda,\mu}(z) = (-1)^n
M_{\lambda,\mu-n\lambda}(z)\;.
\end{align}
Hence, \eqref{E:G-dx} follows.
This completes the proof of Lemma \ref{L:Green}.
\end{proof}

\begin{remark}
Note that in general, the function $x\mapsto G_\beta(t,x)$ is not differentiable
at $x=0$. But we have
\[
\frac{\partial^n}{\partial x^n} G_\beta(t,0-) = (-1)^n
\frac{\partial^n}{\partial x^n} G_\beta(t,0+)
=
\frac{(-1)^n t^{\Ceil{\beta}-1-(n+1)\beta/2}}{2}
\Gamma\left(\Ceil{\beta}-\frac{\beta(n+1)}{2}\right)^{-1}\;,
\]
because $M_{\beta/2,\Ceil{\beta}-n\beta/2}(0-)= M_{\beta/2,\Ceil{\beta}-n\beta/2}(0+) =
\Gamma\left(\Ceil{\beta}-\frac{\beta(n+1)}{2}\right)^{-1}$.
When $\beta=1$ and $n\ge 1$ is an odd integer,
then
$M_{1/2,1-n/2}(0) = \Gamma\left(\frac{1-n}{2}\right)^{-1}=0$, which explains
why the heat kernel function \eqref{E:HeatG} is
smooth at $x=0$.
\end{remark}

\begin{figure}[h!tbp]
\centering
\subfloat[$\beta=6/5$.]{
\includegraphics[scale=0.4]{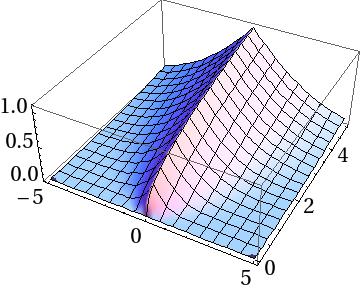}}
\subfloat[$\beta=3/2$.]{
\includegraphics[scale=0.4]{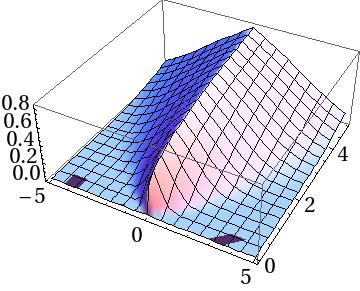}}
\subfloat[$\beta=15/8$.]{
\includegraphics[scale=0.4]{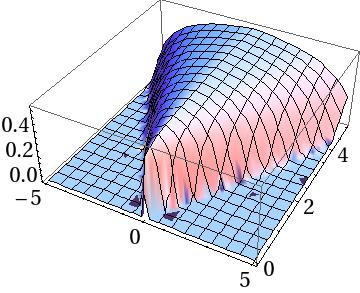}}
\caption{Graphs of the Green functions $G_\beta(t,x)$ for $1<\beta<2$.}
\label{F6:GreenST-FD}
\end{figure}

\begin{remark}[Wave equation case $\beta=2$]\label{R:Wave}
By definition of $M_{\lambda,\mu}(z)$ in \eqref{E:2p-Mainardi}, the parameter $\lambda$ should be strictly less than $1$. Hence, the Green functions $G_\beta(t,x)$ and $G_\beta^*(t,x)$ in
\eqref{E:Green} do not cover the case where $\beta=2$.
However, the wave equation case $\beta=2$ does be a limiting case as $\beta\uparrow 2$, which can be seen from Figure \ref{F6:GreenST-FD}. Another way to see this is through the Fourier transform
\eqref{E:G-Fourier}.
By letting $\beta=2$ in \eqref{E:G-Fourier}, one has that
\begin{align*}
 \int_\R \ud x \; e^{-i\xi x} G_2(t,x) &= t E_{2,2}(-t^2\xi^2) = \frac{\sin(t\xi)}{\xi},\\
\int_\R \ud x \; e^{-i\xi x} G_2^*(t,x) &=E_{2,1}(-t^2\xi^2)=\cos(t\xi),
\end{align*}
which equal the Fourier transforms of the wave kernel functions: $\frac{1}{2}\Indt{|x|\le t}$ and $\frac{1}{2}\left(\delta_{t}(x)+\delta_{-t}(x)\right)$ respectively. Hence, in the limiting case, we
have \eqref{E:WaveG}.
\end{remark}

\begin{figure}[htbp]
\centering
\subfloat[Graphs of $G_\beta(1,x) $ in the linear scale.]{
\includegraphics[scale=0.62]{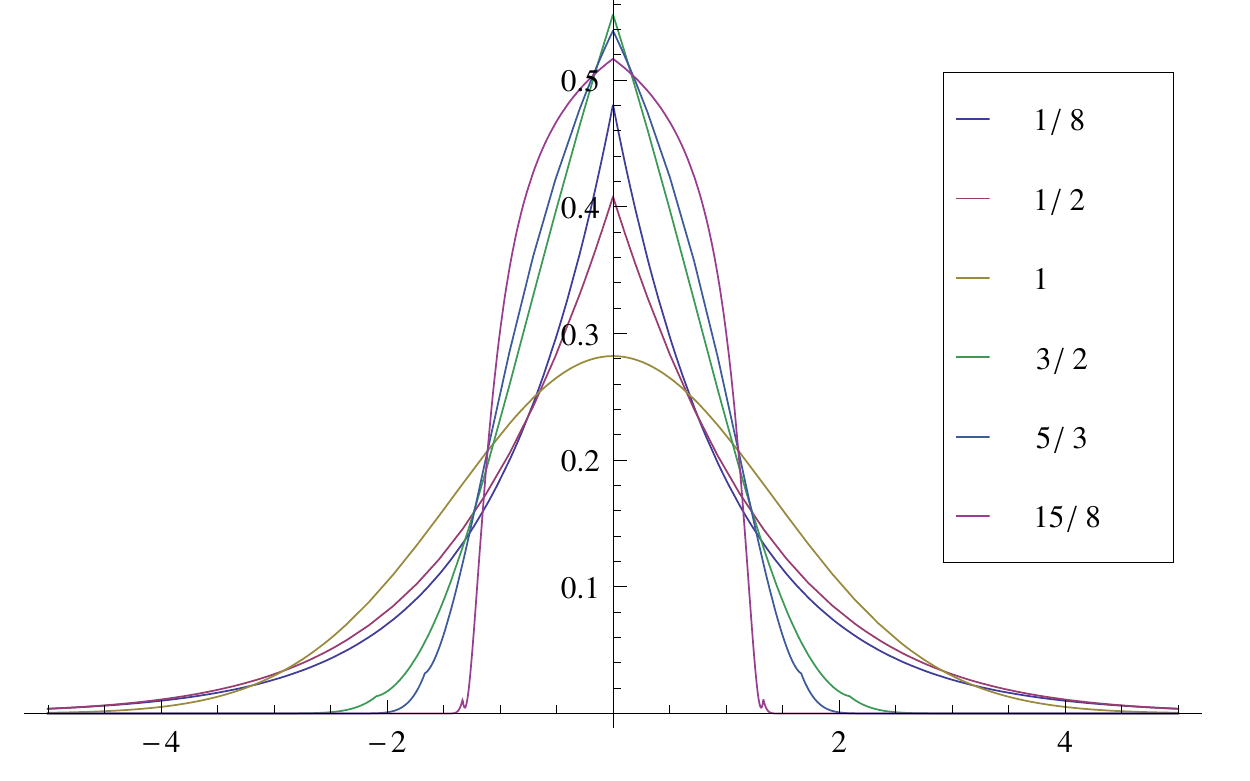}}
\subfloat[Graphs of $\log_{10} G_\beta(1,x)$.]{
\includegraphics[scale=0.67]{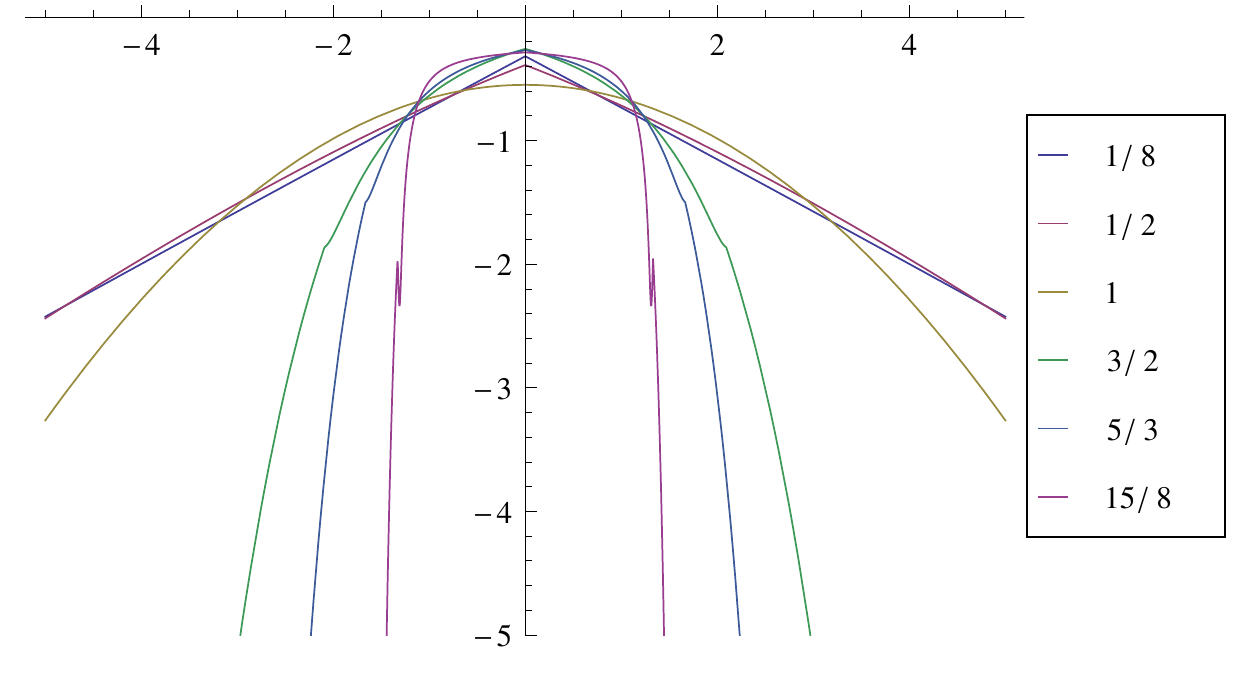}}
\caption{Some graphs with $\beta=1/8$, $1/2$,
$1$, $3/2$, $5/3$ and $15/8$.}
\label{F6:TimeFrac-Green}
\end{figure}

We draw some of these Green functions $G_\beta(1,x)$ in Figure \ref{F6:TimeFrac-Green}.
The range of $x$ is from  $-5$ to $5$.
From these graphs, one can see that when $\beta$ tends to $2$, the Green
function tends to the wave kernel function $\frac{1}{2}\Indt{|x|\le 1}$.
Note that these graphs are plotted by concatenating the truncated summations
for $n \le 23$ in the asymptotic representation \eqref{E:AsymG}, and
hence there are some truncation errors, which can be seen, in these graphs.

In Figure \ref{F6:GreenST-FD}, we draw some Green functions in space-time
coordinates for the fast diffusion equations ($\beta\in\:]1,2[\:$). The ranges for $x$ and
$t$ are $[-5,5]$ and $]0,5]$, respectively. When
$\beta$ tends to $2$, these graphs become closer to the wave kernel function $G_2(t,x)=\frac{1}{2}\Indt{|x|\le t}$.
%

At the end of this section, we list some technical results used in the proof of Lemma \ref{L:Green}.
\begin{lemma}\label{L:Mf-Mom}
The following integral holds:
\[
\int_0^\infty \ud x\: x^a M_{\lambda,\mu}(x) =
\frac{\Gamma\left(a+1\right)}{\Gamma\left(\lambda a + \mu\right)}\;,
\qquad \text{for $a>-1$, $\lambda\in [0,1[\;$, and $\mu\in\bbC$\:.}
\]
\end{lemma}
\begin{proof}
By the integral representation of the Wright funct
ion,
\[
W_{\lambda,\mu}(z) = \frac{1}{2\pi i}\int_{H_a}\frac{\ud t}{t^{\mu}}\:
\exp\left(t+z t^{-\lambda}\right) \;,\qquad \lambda>-1,\quad
\mu\in\bbC\:,
\]
where $H_a$ denotes the Hankel contour
(see \cite[(F.2), on p. 238]{Mainardi10Book} for
more details).
Notice that $M_{\lambda,\mu}(x)=W_{-\lambda,\mu-\lambda}(-x)$.
Then
\begin{align*}
\int_0^\infty \ud x\: x^a M_{\lambda,\mu}(x)
&= \int_0^\infty \ud x\: x^a \left[\frac{1}{2\pi
i}\int_{H_a}\frac{\ud
t}{t^{\mu-\lambda}}e^{t-x t^{\lambda}}\right]=
\frac{1}{2\pi i}\int_{H_a} \ud t\: \left[\int_0^\infty \ud x\:
e^{-x t^\lambda} x^a
\right] \frac{e^t}{t^{\mu-\lambda}}\\
&=
\frac{1}{2\pi i} \int_{H_a} \ud t\:
\frac{\Gamma( a+1)}{t^{\lambda a +\lambda}}
\frac{e^t}{t^{\mu-\lambda}}
=\frac{\Gamma(a+1)}{2\pi i} \int_{H_a}\ud t\:
\frac{e^t}{t^{\lambda a+\mu}}
=
\frac{\Gamma(a+1)}{\Gamma(\lambda a +\mu)}\;,
\end{align*}
where we have used the definition of the Gamma
function in the third equality
(which requires that $a>-1$) and in the last step
we have used the Hankel
integral representation of the Gamma function
$\frac{1}{\Gamma(z)} =\frac{1}{2\pi i} \int_{H_a}\ud t\: e^t t^{-z}$;
see, e.g., \cite[5.9.1, on p. 139]{NIST2010}.
\end{proof}

\begin{lemma}\label{L:Mf-Fourier}
The Fourier transform of the function $M_{\lambda,
\mu}(|x|)$ is
\[
\calF\left[\frac{1}{2} M_{\lambda,\mu}(|\cdot|)\right](\xi)
=\int_\R\ud x\: e^{-i x \xi} M_{\lambda,\mu}(|x|)
=E_{2\lambda,\mu}\left(-\xi^2\right),
\quad \text{for all $\lambda\in [0,1[\;$ and $\mu\in\bbC$\:.}
\]
\end{lemma}
\begin{proof}
By developing in series the cosine function and
the moment formula in Lemma
\ref{L:Mf-Mom},
\begin{align*}
\calF\left[\frac{1}{2}M_{\lambda,\mu}(|\cdot|)\right](\xi)
&= \int_0^\infty \ud x\: \cos(\xi x) M_{\lambda,\mu}(x) =\sum_{n=0}^\infty
(-1)^n \frac{\xi^{2n}}{(2n)!} \int_0^\infty\ud x\: x^{2n}
M_{\lambda,\mu}(x)\\
&=\sum_{n=0}^\infty \frac{\left(-\xi^2\right)^n}{\Gamma\left(2n\lambda +
\mu\right)}= E_{2\lambda,\mu}\left(-\xi^2\right).
\end{align*}
\end{proof}
\begin{theorem}[Bernstein's theorem {\cite[Theorem 12a]{Widder41LaplaceTr}}]
\label{T:Bernstein}
A necessary and sufficient condition that $f(x)$
should be completely monotonic
in $0\le x<+\infty$ is that $f(x) = \int_0^\infty
e^{-x t}\ud \alpha(t)$,
where $\alpha(t)$ is bounded and non-decreasing
and the integral converges
for $0\le x<+\infty$.
\end{theorem}

\section{Calculations of $\calK(t,x)$ and proof of Theorem \ref{T:K-bounds}}
\label{S:K-bounds}
Let $G: \R_+\times\R^d\mapsto \R$ with $d\in\bbN$, $d \ge 1$ be a Borel measurable function.

\begin{assumption}\label{A6:FD}
The function $G: \R_+\times\R^d\mapsto \R$ has the following properties:
\begin{enumerate}[(1)]
 \item There is a nonnegative function $\calG(t,x)$, called {\it reference kernel function}, and constants $C_0>0$,
$\sigma<1$ such that
\begin{align}\label{E:BdScSG}
G(t,x)^2\le \frac{C_0}{t^\sigma}\; \calG(t,x)\;,\quad\text{for all $(t,x)\in\R_+\times\R^d$.}
\end{align}
\item The reference kernel function $\calG(t,x)$ satisfies the following {\it sub-semigroup property}: for some constant $C_1>0$,
\begin{align}\label{E:SG}
\int_{\R^d} \ud y\: \calG\left(t,x-y\right)\calG\left(s,y\right)\le
C_1\; \calG\left(t+s,x\right)\;,\quad\text{for all $t, s>0$ and
$x\in\R^d$.}
\end{align}
\end{enumerate}
\end{assumption}
Define
\[
\calL_0\left(t,x\right) := G(t,x)^2\;,
\quad \text{for all $(t,x)\in\R_+^*\times \R^d$}\;.
\]
Recall that ``$\star$'' denotes the convolution in both space and time variables (space-time convolution). For all $n\in\bbN^*:=\bbN\setminus\{0\}$
and all $(t,x)\in\R_+^*\times\R^d$, define
\begin{align}\notag
\calL_n\left(t,x\right)&:=
\underbrace{\left(\calL_0\star \cdots\star\calL_0\right)}_{\text{$n+1$
times of $\calL_0$}}(t,x),\\
\calK\left(t,x\right)&:= \sum_{n=0}^\infty
\calL_n\left(t,x\right).
\label{E:K-A}
\end{align}
Denote
\[
B_n\left(t;\sigma,C_0,C_1\right):=
C_0^{n} C_1^{n-1}
\frac{\Gamma\left(1-\sigma\right)^{n}}
{\Gamma\left(n(1-\sigma)\right) }\: t^{n(1-\sigma)
-1}\;,\qquad
\text{for all $n\ge 0$}\:.
\]
For simplicity, we write $B_n\left(t;\sigma,C_0,C_1\right)$ simply by
$B_n(t)$.

\begin{proposition}\label{P:K-A}
Under Assumption \ref{A6:FD},  the following properties are true:\\
(i) $\calL_n(t,x)$ is nonnegative and satisfies the following inequality
\begin{align}\label{E:K-i-A}
 \calL_n(t,x) \le \:B_{n+1}(t) \calG(t,x)
\;,\quad\text{for all $n\ge 0$ and $(t,x)\in \R_+^*\times\R^d$.}
\end{align}
Moreover, \eqref{E:K-i-A} becomes an equality if both \eqref{E:BdScSG} and
\eqref{E:SG} are equalities.\\
(ii) For all $t>0$ and $\lambda>0$, the following
 series $\sum_{n=1}^\infty
\calL_n(t,x) $ converges uniformly over
$x\in\R^d$ and hence $\calK(t,x)$ in \eqref{E:K-A}
 is well defined.\\
(iii) $B_n(t)$ are nonnegative and for all $m\in\bbN^*$,  $\sum_{n=0}^\infty
B(t)^{1/m}<+\infty$. \\
(iv) For all $t\ge 0$ and $x\in\R^d$,
\begin{align}\label{E:UpBd-K-Mit-A}
\calK\left(t,x\right)& \le \calG(t,x)\;
\frac{\gamma}{t^{\sigma}}E_{1-\sigma,1-\sigma}
\left(
\gamma t^{1-\sigma}
\right)\\
&\le
\frac{C }{t^{\sigma}} \;\calG(t,x)\;
\left(
1+t^{\sigma} \exp\left(\gamma^{\frac{1}{1-\sigma}}
\; t\right)
\right)\;,
\label{E:UpBd-K-A}
\end{align}
where $\gamma=C_0 C_1\Gamma(1-\sigma)$ and the
constant
$C=C\left(\sigma,\gamma\right)$ can be chosen as
\begin{align}\label{E:Cst-UpK-A}
C\left(\sigma,\gamma\right):= \gamma \: \sup_{t\ge 0}
\frac{E_{1-\sigma,1-\sigma}\left(\gamma\:
t^{1-\sigma}\right)}{1+t^{\sigma}\exp\left(\gamma^
{\frac{1}{1-\sigma}}
\; t\right)}<+\infty\;.
\end{align}
Moreover, \eqref{E:UpBd-K-Mit-A} becomes equality
if both \eqref{E:BdScSG} and
\eqref{E:SG} are equalities. \\
(v) If there exist a kernel function $\ubar{\calG}(t,x)$ and some constants
$\ubar{C}_0>0$, $\ubar{C}_1>0$, and
$\ubar{\sigma}<1$ such that for all $t, s>0$
and $x\in\R^d$,
\[
G(t,x)^2\ge \ubar{C}_0\;t^{-\ubar{\sigma}}\:
\ubar{\calG}(t,x),
\]
and $\ubar{\calG}(t,x)$ satisfies the {\it sup-semigroup property}
\[
\int_{\R^d} \ud y\:
\ubar{\calG}\left(t,x-y\right)\ubar{\calG}\left(s,
y\right)\ge
\ubar{C}_1\; \ubar{\calG}\left(t+s,x\right),
\]
then for all $t\ge 0$ and $x\in\R^d$,
\begin{align}\label{E:LowBd-K-Mit-A}
\calK\left(t,x\right)& \ge \ubar{\calG}(t,x)\;
\frac{\ubar{\gamma}}{t^{\ubar{\sigma}}}E_{1-\ubar{
\sigma},1-\ubar{\sigma}}
\left(\;
\ubar{\gamma} t^{1-\ubar{\sigma}}
\right)\\
&\ge
\ubar{C}\;\ubar{\calG}(t,x)\;
 \exp\left(\:\ubar{\gamma}^{\frac{1}{1-\ubar{\sigma}}}\; t\right)\;,
\label{E:LowBd-K-A}
\end{align}
where
$\ubar{\gamma}=\ubar{C}_0\ubar{C}_1\Gamma(1-\ubar{
\sigma})$ and
\[
\ubar{C}=\ubar{C}\left(\ubar{\sigma},\ubar{\gamma}
\right):=\ubar{\gamma}\inf_{
t\ge 0}
\frac{E_{1-\ubar{\sigma} , 1-\ubar{ \sigma } }
\left(\ubar{\gamma}\:
t^{1-\ubar{\sigma}}\right)}{t^{\ubar{\sigma}}
\exp\left(\ubar{\gamma}^{\frac{1}{1-\ubar{\sigma}}
}\:
t\right)}>0\:.
\]
\end{proposition}
\begin{proof}
(i) The non-negativity is clear. The case $n=0$
is trivially true.
Suppose that the relation \eqref{E:K-i-A} holds
up to $n-1$. Then by
the Beta integral,
\begin{align*}
\calL_n(t,x)
&= \int_0^t\ud s \int_{\R^d}\ud y\: \calL_{n-1}\left(t-s,x-y\right)
G^2\left(s,y\right)\\
&\le
C_0 \int_0^t \ud s \: B_n(t-s) s^{-\sigma} \int_{\R^d}\ud y\:
\calG\left(t-s,x-y\right)\calG\left(s,y\right)\\
&\le
C_0^{n+1} C_1^{n} \;\calG(t,x)
\frac{\Gamma(1-\sigma)^n}{\Gamma(n(1-\sigma))}\int_0^t\ud s\:
(t-s)^{n(1-\sigma)-1}s^{-\sigma}\\
&=
C_0^{n+1} C_1^{n}
\frac{\Gamma(1-\sigma)^{n+1}}{\Gamma((n+1)(1-\sigma))}\;
t^{(n+1)(1-\sigma)-1}\;\calG(t,x)\\
&= B_{n+1}(t) \calG(t,x)\;.
\end{align*}

(ii) It is a special case of (iii). (iii) The non-
negativity is clear.
By \eqref{E:K-i-A},
\[
\calL_n(t,x)\le B_{n+1}(t)
t^{-\sigma}\sup_{x\in\R^d}\calG(t,x)<+\infty\;.
\]
Thus, if the series  $\sum_n B_n(t)^{1/m}$
converges, then it does so
uniformly over $x\in\R^d$.
Denote $\beta :=1-\sigma$.
Use the ratio test
\[
\left(\frac{B_n(t)}{B_{n-1}(t)}\right)^{1/m} =
\left(C_0 C_1
\Gamma\left(\beta\right) t^\beta\right)^{1/m}
\left(
\frac{\Gamma\left((n-1)(1-\sigma)\right)}
{\Gamma(n(1-\sigma))}
\right)^{1/m} \;.
\]
By the asymptotic expansion of the Gamma function
(\cite[5.11.2, on p.
140]{NIST2010}),
\[
\frac{\Gamma\left((n-1)(1-\sigma)\right)}{\Gamma\left(n(1-\sigma) \right) }
\approx \left(\frac{e}{\beta}\right)^\beta
\left(1-\frac{1}{n}\right)^{(n-1)\beta} \frac{1}{n^\beta}
\approx \frac{1}{(\beta n)^\beta}
\]
for large $n$. Now clearly, $\beta >0$ since $\sigma<1$. Hence, for all $t>0$
and for large $n$,
\[
\left(
\frac{B_n(t)}
{B_{n-1}(t)}
\right)^{1/m}
\approx
\left(
C_0 C_1
\Gamma\left(\beta\right) t^\beta\right)^{1/m} \frac{1}{(\beta n)^{\beta/m}},
\]
which tends to zero as $n\rightarrow+\infty$.

(iv) The bound \eqref{E:UpBd-K-Mit-A} is because
\begin{align}\label{E:Efrom1}
\sum_{k=1}^{\infty} \frac{z^k}{\Gamma(\alpha k)}
= z E_{\alpha,\alpha}(z)\;,
\end{align}
and the bounds in (i):
\begin{align*}
\calK\left(t,x\right)
&\le  \calG(t,x) \sum_{n=1}^\infty B_{n}\left(t\right)
=
\frac{1}{C_1 t}\calG(t,x) \sum_{n=1}^\infty \frac{\left(C_0C_1
\Gamma(1-\sigma)\;t^{1-\sigma}\right)^n}{\Gamma(n(1-\sigma))}\\
&=
C_0 \;\Gamma(1-\sigma)\; t^{-\sigma}\;
\calG(t,x)\;E_{1-\sigma,1-\sigma}\left(C_0C_1 \Gamma(1-\sigma)
t^{1-\sigma}\right).
\end{align*}
As for \eqref{E:UpBd-K-A}, we only need to show
that the constant $C$ defined
in \eqref{E:Cst-UpK-A} is finite.
Let
\[
f(t) = \frac{E_{1-\sigma,1-\sigma}\left(\gamma\:
t^{1-\sigma}\right)}{1+t^{\sigma}\exp\left(\gamma^
{\frac{1}{1-\sigma}}\;
t\right)}\;.
\]
By Lemma \ref{L5:Eab} with real nonnegative
value $z= \gamma t^{1-\sigma}$ and $p=1$:
\[
\gamma E_{1-\sigma,1-\sigma}\left(\gamma\; t^{1-\sigma}\right)=
\frac{1}{1-\sigma} \gamma^{\frac{1}{1-\sigma}}\;
t^{\sigma}
\exp\left(\gamma^{\frac{1}{1-\sigma}}\; t\right) +
O\left(\frac{1}{|t|^{1-\sigma}}\right)
\:,\qquad t\rightarrow +\infty\:,
\]
we see that $\lim_{t\rightarrow+\infty} f(t) = \frac{1}{1-\sigma}
\gamma^{\frac{1}{1-\sigma}}$.
Then because the Mittag-Leffler function is an
entire
function on complex plain \cite[Theorem 4.1, p.
68]{Diethelm10AFDE},
we can conclude that $\sup_{t\ge 0} f(t)<+\infty$.

(v) The proof is similar to (i) and (iv). We only
need to show that
$\ubar{C}$ is strictly positive. Because the function
$g(t) =E_{1-\ubar{\sigma},1-\ubar{\sigma}}^{-1}\left(\ubar{\gamma}\:
t^{1-\ubar{\sigma}}\right)\:
t^{\ubar{\sigma}}\exp\left(\ubar{\gamma}^{\frac{1}
{1-\ubar{\sigma}}}\;
t\right)$ is continuous over $t\in [0,+\infty]$
with $g(0)=0$ and
$\lim_{t\rightarrow +\infty} g(t)=(1-\ubar{\sigma})\:
\ubar{\gamma}^{\frac{\ubar{\sigma}}{\ubar{\sigma}-1}}<+\infty$,
this function is bounded from above for $t\in [0,+\infty]$ and hence
$\inf_{t\ge 0} g^{-1}(t)>0$.
This completes the proof of Proposition \ref{P:K-A}.
\end{proof}

\begin{example} \label{Ex:SHE}
For the heat kernel
$p_\nu(t,x)=(2\pi\nu t)^{-1/2}\exp\left(-\frac{
x^2}{2\nu t}\right)$ with
$\nu>0$,
Assumption \ref{A6:FD} holds with both inequalities \eqref{E:BdScSG} and
\eqref{E:SG} replaced by equalities, and
\[
C_0=\frac{1}{\sqrt{4\pi\nu}}, \quad \sigma=\frac{
1}{2},\quad \calG(t,x)=
p_{\nu/2}(t,x),\quad C_1=1.
\]
Then, $\gamma=(4\nu)^{-1/2}$. Therefore, by \eqref
{E:ML-HH} and
$\Erfc(-x)=2\Phi(\sqrt{2} x)$ where $\Phi(x)$ is
the distribution function of
the standard normal distribution, Proposition \ref{P:K-A} implies that
\begin{align*}
\calK(t,x)&=  \frac{\calG(t,x)}{\sqrt{4\nu t}}
\left[
\frac{1}{\sqrt{\pi}}+ \frac{\sqrt{t}}{\sqrt{4\nu}}
\Erfc\left(-\frac{\sqrt{t}}{\sqrt{4\nu}}\right)
e^{\frac{t}{4\nu}}
\right]=
p_{\frac{\nu}{2}}(t,x)
\left[
\frac{1}{\sqrt{4\pi\nu t}}+ \frac{1}{2\nu}
\Phi\left(\frac{\sqrt{t}}{\sqrt{2\nu}}\right) e^{\frac{t}{4\nu}}
\right],
\end{align*}
which recovers the results in \cite{ChenDalang13Heat}.
\end{example}

\begin{example}\label{Ex:DoubleSHE}
Let us consider the following SPDE
\[
\begin{cases}
\left(\frac{\partial}{\partial t}-\frac{\partial^2
}{\partial x^2}\right)^2
u(t,x) = \rho(u(t,x))\W(t,x),&(t,x)\in\R_+^*\times
\R, \\
u(0,\cdot) = \mu(\cdot).
\end{cases}
\]
The Green function is
$G(t,x)= \frac{\sqrt{t}}{\sqrt{4\pi}} \exp\left(-\frac{x^2}{4t}\right)$; see
\cite[Section 9.2.5-2]{Polyanin02LPDE}.
Assumption \ref{A6:FD} holds with both inequalities \eqref{E:BdScSG} and
\eqref{E:SG} replaced by equalities, and
\[
C_0=\frac{1}{\sqrt{8\pi}}, \quad \sigma=-\frac{3}{2},\quad \calG(t,x)=
\frac{1}{\sqrt{2\pi t}}\exp\left(-\frac{x^2}{2 t}\right),\quad C_1=1.
\]
Then, $\gamma=\frac{1}{\sqrt{8\pi}} \Gamma(5/2) =
\frac{3\sqrt{2}}{16}$.
Therefore, Proposition \ref{P:K-A} implies that
\[
\calK(t,x)= \frac{3\sqrt{2}}{16}\: t^{3/2}\:
\calG(t,x)\;E_{5/2,5/2}\left(\frac{3\sqrt{2}}{26}\: t^{5/2} \right).
\]
In particular, if $\rho(u)=u$, then $\E\left(u^2(
t,x)\right)= J_0^2(t,x)+
(J_0^2\star \calK)(t,x)$, where $J_0(t,x) = (\mu*
G(t,\cdot))(x)$.
Note that the initial data can be more general
than the SHE \eqref{E:Heat}: It can be any distribution $\mu$ such that it is
the (distributional) derivative
of some measures in $\calM_H(\R)$, i.e, for some
$\mu_0\in\calM_H( \R)$, $\mu=\mu_0'$.
More details of this SPDE, which will not be pursued here, are left
to interested readers.
\end{example}

Here are three natural choices of the reference kernel functions $\calG(t,x)$:
\begin{enumerate}
 \item[(1)] The Gaussian kernel function
\begin{align*}
\calG_g(t,x)&:= \left(4\pi t\right)^{-d/2}\exp\left(-\frac{|x|^2}{4 t}\right),\quad\text{for all $(t,x)\in\R_+\times\R^d$, $d\ge 1$,}
\end{align*}
where $|x|^2=x_1^2+\cdots+x_d^2$;
\item[(2)] The Poisson kernel function:
\[
\calG_p(t,x):=c_n \: \frac{t}{\left(t^2+|x|^2\right)^{(d+1)/2}},\quad\text{for all $(t,x)\in\R_+\times\R^d$, $d\ge 1$,}
\]
where $c_n=\pi^{-(n+1)/2}\: \Gamma((n+1)/2)$;
\item[(3)] The exponential kernel function $\calG_{e,\beta}(t,x)$ defined in \eqref{E:Geb}.
\end{enumerate}

Clearly, we have the following scaling properties for these reference kernel functions:
\begin{align*}
 \calG_g(t,x)&=t^{-d/2}\:\calG_g(1,t^{-1/2} \:x),\\
\calG_p(t,x)&=t^{-d}\: \calG_p(1,t^{-1} \:x),\\
\calG_{e,\beta}(t,x)&=t^{-\beta/2} \calG_{e,\beta}(1,t^{-\beta/2}x).
\end{align*}
Both $\calG_g(t,x)$ and $\calG_p(t,x)$ satisfy part (2) of Assumption \ref{A6:FD} with $C_1 =1$ and ``$\le$'' replaced by ``$=$''.
By Lemma \ref{L:SG-Geb} below, $\calG_{e,\beta}(t,x)$ satisfies part (2) of Assumption \ref{A6:FD} with $C_1 =\widehat{C}_\beta$, where $\widehat{C}_\beta$ is defined in \eqref{E:HatC}.

\begin{proposition}[Gaussian reference kernel]
\label{P:GaussRef}
Suppose the function $G:\R_+\times\R^d\mapsto \R$
satisfies the
following two properties:
\begin{enumerate}[(i)]
 \item The scaling property:  for some constants
$\gamma_1\in\R$ and
$\gamma_2\ge 1/2$,
\[
G(t,x) = t^{\gamma_1} G\left(1,\frac{x}{t^{\gamma_2}}\right),\quad \text{for
all $(t,x)\in\R_+\times\R^d$;}
\]
\item The function $x\mapsto G(1,x)$ is bounded
such that
$\sup_{x\in \R^d} \frac{G(1,x)^2}{\calG_g(1,x)} < +\infty$.
\end{enumerate}
Then $G(t,x)$ satisfies Assumption \ref{A6:FD}
with
$\calG(t,x)=\calG_g(t^{2\gamma_2},x)$ and
\begin{align}\label{E:GR-Const}
C_0=\sup_{x\in \R^d} \frac{G(1,x)^2}{\calG_g(1,x)}
,
\quad
C_1=2^{d\left(\gamma_2 -1/2\right)},
\quad\text{and}\quad
\sigma=- (2\gamma _1+d\gamma_2).
\end{align}
\end{proposition}
\begin{proof}
Notice that by the scaling properties of $\calG_g(t,x)$ and $G(t,x)$, we have that
\begin{align*}
 \sup_{(t,x)\in\R_+\times\R^d}\frac{G(t,x)^2}{t^{-\sigma}\calG(t,x)}=
\sup_{(t,x)\in\R_+\times\R^d}\frac{G(t,x)^2}{t^{-\sigma}\calG_g(t^{2\gamma_2},x)}=
\sup_{y\in\R^d} \frac{G(1,y)^2}{\calG_g(1,y)},
\end{align*}
which is finite by (ii). Hence, part (1) of Assumption \ref{A6:FD} is satisfied
with the constants $C_0$ and
$\sigma$ defined in \eqref{E:GR-Const}.
As for part (2) of Assumption \ref{A6:FD}, by the
semigroup property of $\calG_g(t,x)$, we have
\begin{align*}
\int_{\R^d} \ud y\: \calG(t,x-y)\calG(s,y)
&=
\int_{\R^d}\ud y\:
\calG_g\left(t^{2\gamma_2},y\right)
\calG_g\left(s^{2\gamma_2},x-y\right)=\calG_g\left(t^{2\gamma_2}+s^{2\gamma_2},
x\right).
\end{align*}
Notice that the function $[0,1]\ni r \mapsto f(r)
= (1-r)^{2\gamma_2} +
r^{2\gamma_2}$ is
convex because $2\gamma_2\ge 1$. By solving $f'(r)=0$, we find that
\[
\min_{r\in [0,1]} f(r) = f(1/2) = 2^{1-2\gamma_2}\;,\quad\text{and}\quad
\max_{r\in [0,1]} f(r) = f(1) = f(0) = 1\;.
\]
Hence,
\begin{align}\label{E:t-s+s_F}
 2^{1-2\gamma_2} (t+s)^{2\gamma_2} \le
t^{2\gamma_2} + s^{2\gamma_2} = (t+s)^{2\gamma_2}
f\left(\frac{s}{t+s}\right)
\le (t+s)^{2\gamma_2}
\end{align}
Finally,
\begin{align*}
\calG_g\left(t^{2\gamma_2}+s^{2\gamma_2},x\right)
&=
\left[4\pi \left(t^{2\gamma_2}+s^{2\gamma_2} \right)\right]^{-d/2}
\exp\left(-\frac{|x|^2}{4\left(t^{2\gamma_2}+s^{2\gamma_2} \right)}\right)\\
&\le
2^{d(\gamma_2-1/2)}
\left[4\pi \left((t+s)^{2\gamma_2} \right)\right]^{-d/2}
\exp\left(-\frac{|x|^2}{4 (t+s)^{2\gamma_2}}\right
)\\
&= 2^{d(\gamma_2-1/2)}\; \calG\left((t+s)^{2\gamma_2},x \right),
\end{align*}
which completes the proof of Proposition \ref{P:GaussRef}.
\end{proof}

We will not use the Poisson reference kernel in this paper. We prove the following result for the future reference.
\begin{proposition}[Poisson reference kernel]
\label{P:PoissonRef}
Suppose the function $G:\R_+\times\R^d\mapsto \R$
satisfies the
following two properties:
\begin{enumerate}[(i)]
 \item The scaling property: For some constants $\gamma_1\in\R$ and
$0<\gamma_2\le 1$,
\[
G(t,x) = t^{\gamma_1} G\left(1,\frac{x}{t^{\gamma_2}}\right),\quad \text{for all $(t,x)\in\R_+\times\R^d$;}
\]
\item The function $x\mapsto G(1,x)$ is bounded
such that
$\sup_{x\in \R^d} \frac{G(1,x)^2}{\calG_p(1,x)} < +\infty$.
\end{enumerate}
Then $G(t,x)$ satisfies Assumption \ref{A6:FD}
with
$\calG(t,x)=\calG_g(t^{\gamma_2},x)$ and
\begin{align}
\label{E:PR-Const}
C_0=\sup_{x\in \R^d} \frac{G(1,x)^2}{\calG_p(1,x)}
,
\quad
C_1=2^{1-\gamma_2},
\quad\text{and}\quad
\sigma=-(2\gamma_1+d \gamma_2 ).
\end{align}
\end{proposition}
\begin{proof}
Notice that by the scaling properties of $\calG_g(t,x)$ and $G(t,x)$, we have
that
\begin{align*}
 \sup_{(t,x)\in\R_+\times\R^d}\frac{G(t,x)^2}{t^{-\sigma}\calG(t,x)}=
\sup_{(t,x)\in\R_+\times\R^d}\frac{G(t,x)^2}{t^{-\sigma}\calG_p(t^{\gamma_2},x)}=
\sup_{y\in\R^d} \frac{G(1,y)^2}{\calG_p(1,y)},
\end{align*}
which is finite due to (ii). Hence, part
(1) of Assumption \ref{A6:FD} is satisfied with
the constants $C_0$ and
$\sigma$ defined in \eqref{E:PR-Const}.
As for part (2) of Assumption \ref{A6:FD}, by the
semigroup property of
$\calG_g(t,x)$, we have
\begin{align*}
\int_{\R^d} \ud y\: \calG(t,x-y)\calG(s,y)
&=
\int_{\R^d}\ud y\:
\calG_p\left(t^{2\gamma_2},y\right)
\calG_p\left(s^{2\gamma_2},x-y\right)=\calG_p\left(t^{2\gamma_2}+s^{2\gamma_2},
x\right).
\end{align*}
Then because $0<\gamma_2\le 1$,
\begin{align}\label{E:t-s+s_S}
 (t+s)^{\gamma_2}\le t^{\gamma_2}+s^{\gamma_2} \le
2^{1-\gamma_2} (t+s)^{\gamma_2}.
\end{align}
Therefore,
\begin{align*}
\calG_p\left(t^{\gamma_2}+s^{\gamma_2},x\right)
&= \frac{c_n
\left(t^{\gamma_2}+s^{\gamma_2}\right)}{\left((t^{\gamma_2}+s^{\gamma_2})^2+|x|^2
\right)^{(d+1)/2}}\\
&\le
\frac{c_n
\:2^{1-\gamma_2}\:\left(t+s\right)^{\gamma_2}}{\left((t+s)^{2\gamma_2}+|x|^2
\right)^{(d+1)/2}} = 2^{1-\gamma_2}\: \calG_p((t+s)^{\gamma_2},x),
\end{align*}
which completes the proof of Proposition \ref{P:PoissonRef}.
\end{proof}

\begin{proposition}[Exponential reference kernel]
\label{P:ExpRef}
Let $\beta\in\;]0,2]$.
Suppose the function $G:\R_+\times\R\mapsto \R$
satisfies the
following two properties:
\begin{enumerate}[(i)]
 \item The scaling property:
\[
G(t,x) = t^{-\beta/2} G\left(1,\frac{x}{t^{\beta/2}}\right),\quad \text{for
all $(t,x)\in\R_+\times\R$;}
\]
\item The function $x\mapsto G(1,x)$ is bounded
such that
$\sup_{x\in \R} \frac{G(1,x)^2}{\calG_{e,\beta}(1,x)} < +\infty$.
\end{enumerate}
Then $G(t,x)$ satisfies Assumption \ref{A6:FD}
with
$\calG(t,x)=\calG_{e,\beta}(t,x)$ and
\begin{align}\label{E:ER-Const}
C_0=\sup_{x\in \R} \frac{G(1,x)^2}{\calG_{e,\beta}(1,x)}
,
\quad
C_1=\widehat{C}_\beta,
\quad\text{and}\quad
\sigma=\frac{\beta}{2},
\end{align}
where $\widehat{C}_\beta$ is defined in \eqref{E:HatC}.
\end{proposition}

\begin{proof}
Notice that by the scaling properties of $\calG_{e,\beta}(t,x)$ and $G(t,x)$, we have
that
\begin{align*}
 \sup_{(t,x)\in\R_+\times\R}\frac{G(t,x)^2}{t^{-\sigma}\calG(t,x)}=
\sup_{y\in\R} \frac{G(1,y)^2}{\calG_{e,\beta}(1,y)},
\end{align*}
which is finite due to (ii). Hence, part
(1) of Assumption \ref{A6:FD} is satisfied with
the constants $C_0$ and
$\sigma$ defined in \eqref{E:ER-Const}.
Part (2) of Assumption \ref{A6:FD} is due to Lemma \ref{L:SG-Geb} below with $C_1=\widehat{C}_\beta$.
This completes the proof of Proposition \ref{P:ExpRef}.
\end{proof}

Now we apply Proposition \ref{P:K-A} to the Green functions $G_\beta(t,x)$ with $0<\beta<2$.
More precisely, we will apply Proposition \ref{P:GaussRef} (resp. \ref{P:ExpRef}) with $\G_\beta(t,x)$ defined in \eqref{E:calG} in the case of fast (resp. slow) diffusions for the upper bounds
of $\calK(t,x)$, and Proposition \ref{P:GaussRef} with $\ubar{\calG}_\beta(t,x)$ defined in
\eqref{E:calG-L} in the case of slow diffusion for the lower bound of $\calK(t,x)$.
Recall the constants $\Psi_\beta$ and $\ubar{\Psi}_\beta$ defined in \eqref{E_:G-Bdd}
and \eqref{E_:G-Bdd-L}, respectively, and the constant $\widetilde{C}_\beta$ defined in \eqref{E:TildeC}.

\begin{proposition}\label{P:ST-Con}
(1) Proposition \ref{P:K-A} (i) -- (iv) hold for
 $G_\beta(t,x)$ with $\beta\in \: ]0,2[\;$ and
\begin{gather*}
d=1\:,\quad
\sigma=\beta/2 +2(1-\Ceil{\beta}),
\quad
\calG(t,x)=\calG_\beta(t,x)\;, \quad
C_0=\Psi_\beta
\;,\quad
C_1=\widetilde{C}_\beta\;.
\end{gather*}
(2) Proposition \ref{P:K-A} (v) holds for $G_\beta(t,x)$ with $\beta\in \:]0,1[\;$ and
\[
d=1, \quad \ubar{\sigma}= \frac{\beta}{2}-1,\quad
\ubar{\calG}(t,x)=\ubar{\calG}_\beta(t,x),\quad
\ubar{C}_0=\ubar{\Psi}_\beta,\quad
\ubar{C}_1 = 2^{\frac{\beta-2}{4}}.
\]
\end{proposition}
\begin{proof}
(1) We begin with the case where $\beta\in \:]1,2[\:$.
By \eqref{E:Scaling}, $G_\beta(t,x)$ satisfies the scaling property with
$\gamma_2=\beta/2 \ge 1/2$ and $\gamma_1=1-\beta/2$.
Notice that
\[
\Psi_\beta =
\sup_{y\ge 0}
\frac{\sqrt{\pi}}{2}\; \exp\left(y^2/4\right)\;
M_{\beta/2}^2\left(y\right).
\]
Because the parameter $c$ in \eqref{E:AsymG} is
strictly bigger than $2$ (see also Figure \ref{F6:Asy-Parameters}), we see that
\[
\lim_{y\rightarrow+\infty}
\frac{\sqrt{\pi}}{2}\; \exp\left(y^2/4\right)
M_{\beta/2}^2\left(y\right)=0\;.
\]
Since the function $y\mapsto\exp\left(y^2/2\right) M_{\beta/2}^2\left(y\right)$ is an entire function, we
see that the above supremum does exist.
Therefore, one can apply Proposition \ref{P:GaussRef} with $d=1$, $C_0=\Psi_\beta$, $C_1=2^{(\beta-1)/2}$, $\sigma=\beta/2-2<-1$,
and the above $\gamma_1$ and $\gamma_2$.

The proof for the slow diffusion equations can be proved
similarly using Proposition
\ref{P:ExpRef} with
$\gamma_2=-\gamma_1= \beta/2\le 1/2$,
$C_0=\Psi_\beta=\sup_{y\ge 0}\; 1/2\:  e^{y}M_{\beta/2}^2\left(y\right)$,
$C_1=\widehat{C}_\beta$, and $\sigma=\beta/2<1$.

(2) We claim that if $\beta \in\:]0,1[\:$, then for all $(t,x)\in\R_+\times\R$ and $s\ge 0$, we have that
\begin{gather}\label{E:LowSD1}
G_\beta^2(t,x) \ge \ubar{\Psi}_\beta \: t^{1-\beta/2}\:
\ubar{\calG}_\beta(t,x),\\
\label{E:LowSD2}
\int_\R \ud y \: \ubar{\calG}_\beta(t,x-y) \ubar{\calG}_\beta(s,y) \ge
2^{\frac{\beta-2}{4}} \: \ubar{\calG}_\beta(t+s,x).
\end{gather}
By the scaling property \eqref{E:Scaling}, $\sup_{(t,x)\in\R_+\times\R}
\frac{t^{1-\beta/2}\: \ubar{\calG}_\beta(t,x)}{G_\beta^2(t,x)} = \sup_{y\in\R}
\frac{\ubar{\calG}_\beta(1,y)}{G_\beta^2(1,y)}$,
which is finite by the same reasoning as above, where the parameter $c$ in \eqref{E:Asym-c}
is strictly less than $2$ in this case.
Thus, $\ubar{C}_0=\ubar{\Phi}_\beta >0$ and \eqref{E:LowSD1}
follows with $\ubar{\sigma}= \beta/2-1<-1/2$.
The inequality \eqref{E:LowSD2} is proved by the semigroup
property of the heat kernel function and \eqref{E:t-s+s_S}. So $\ubar{C}_1=2^{(\beta-2)/4}$.
Then apply Proposition \ref{P:GaussRef}.
This completes the proof of
Proposition \ref{P:ST-Con}.
\end{proof}

At the end of this section, we list two technical results that are used in this section.

\begin{lemma}[Theorem 1.3, p. 32 in \cite{Podlubny99FDE}]\label{L5:Eab}
If $0<\alpha<2$, $\beta$ is an arbitrary complex
number and $\mu$ is an
arbitrary real number such that $\pi\alpha/2<\mu<\pi \wedge (\pi\alpha)$,
then for an arbitrary integer $p\ge 1$ the following expression holds:
\[
E_{\alpha,\beta}(z) = \frac{1}{\alpha} z^{(1-\beta)/\alpha}
\exp\left(z^{1/\alpha}\right)
-\sum_{k=1}^p \frac{z^{-k}}{\Gamma(\beta-\alpha k)
} + O\left(|z|^{-1-p}\right)
,\quad |z|\rightarrow\infty,\quad |\arg(z)|\le \mu\:.
\]
\end{lemma}

\begin{lemma}\label{L:SG-Geb}
Suppose $\beta\in \:]0,2]$. The exponential reference kernel function $\calG_{e, \beta}(t,x)$ defined in \eqref{E:Geb} satisfies the sub-semigroup property, i.e., for all $t\ge0$, $s \ge 0$ and
$x\in\R$,
\[
\left(\calG_{e,\beta}(t,\cdot)* \calG_{e,\beta}(s, \cdot)\right) (x) \le
 \widehat{C}_{\beta} \:
\calG_{e,\beta}\left(t+s, x\right),
\]
where the constant $\widehat{C}_\beta$ is defined in \eqref{E:HatC}.
\end{lemma}

\begin{proof}
 Fix $a> b>0$ and let $\theta=\beta/2$.
Because
\[
\int_\R \frac{1}{4a^\theta b^\theta} \exp\left(- \frac{|x-y|}{a^\theta}
-\frac{|y|}{b^\theta}\right)\ud y
=\frac{1}{2(a^\theta+b^\theta)} \frac{a^\theta\; \exp\left(-\frac{|x|}{a^\theta}\right)-b^\theta\;
\exp\left(-\frac{|x|}{b^\theta}\right)}{a^\theta-b^\theta}\;,
\]
we only need to prove that
\begin{align}\label{E6:SD-SG}
\frac{a^\theta\; \exp\left(-\frac{|x|}{a^\theta}\right)-b^\theta\;
\exp\left(-\frac{|x|}{b^\theta}\right)}{a^\theta-b^\theta}\le
\widehat{C}_{2\theta} \exp\left(-\frac{|x|}{(a+b)^\theta}\right),\quad\text{for all
$(t,x)\in\R_+\times\R$.}
\end{align}
By setting $r=b/a$ and $\eta= |x|/a^\theta$, \eqref{E6:SD-SG} is equivalent to
\[
\frac{e^{-\eta}-r^\theta e^{-\eta/r^\theta}}{1-r^\theta}\le \widehat{C}_{2\theta}\;
\exp\left(-\frac{\eta}{(1+r)^\theta}\right)\;,\quad\text{for all $r\in ]0,1[\;$ and $\eta\ge 0$}\;.
\]
Denote
\[
f\left(r,\eta\right) =
\frac{e^{-\eta}-r^\theta e^{-\eta/r^\theta}}{1-r^\theta}
\exp\left(\frac{\eta}{(1+r)^\theta}\right).
\]
Some simple calculations show that
\[
\lim_{r\rightarrow 0_+} f\left(r,\eta\right) =1
\quad\text{and}\quad
\lim_{r\rightarrow 1_-} f\left(r,\eta\right) =(1+\eta)e^{-(1-2^{-\theta})\eta}\;,
\]
and also
\[
\lim_{\eta \rightarrow 0_+} f\left(r,\eta \right)
=1\quad\text{and}\quad
\lim_{\eta \rightarrow +\infty} f\left(r,\eta \right) =0\;,
\]
Fix $r\in \;]0,1[\;$. By solving
\[
\frac{\partial f\left(r,\eta \right)}{\partial \eta}
=
\frac{ e^{\eta  \left(-r^{-\beta }+(1+r)^{-\beta
   }-1\right)} \left(e^{\eta } \left((1+r)^{\beta}-r^{\beta }\right)- \left((1+r)^{\beta }-1\right) e^{\eta/ r^{\beta
   }}\right)}{(1+r)^\beta(1-r^{\beta})} =0,
\]
which has one finite solution
\[
\eta = -\frac{r^\theta}{1-r^\theta} \log\frac{(1+r)^\theta-1}{(1+r)^\theta-r^\theta},
\]
we find the local maximum of the function $\eta\mapsto f\left(r,\eta\right)$. This local maximum is indeed the global maximum.
Hence
\[
f\left(r,\eta\right)
\le h(r):= \frac{(1+r)^{\theta }}{(1+r)^{\theta }-1}\exp\left[\frac{1-\left(\frac{r}{1+r}\right)^{\theta
   }}{1-r^{\theta }}\log\left(1-\frac{1-r^{\theta }}{(1+r)^{\theta
   }-r^{\theta }}\right)\right]\:.
\]
Because $h'(r)\ge 0$ for all $r\in \; ]0,1[\;$ and $\theta\in \; ]0,1]$, we have that
\[
h(r)\le
\lim_{r\rightarrow 1} h(r) = \frac{2^{\theta }}{2^{\theta }-1} \exp\left(-\frac{1}{2^{\theta}}\right)=
\widehat{C}_{2\theta}.
\]
Therefore, $f\left(r,\beta\right)\le \widehat{C}_{2\theta}$.
 This proves \eqref{E6:SD-SG}.

Apply \eqref{E6:SD-SG} with
$a=\left(t\vee s\right)^{\beta/2}$, $b=\left(t\wedge s\right)^{\beta/2}$ and $\theta=\beta/2 \in \: ]0,1]$, and use $(t+s)^{\beta/2}\le t^{\beta/2}+s^{\beta/2}$ to have
\begin{align*}
\int_\R \ud y \: \calG_{e,\beta}(t,x-y)\: \calG_{e,\beta}(s,y)
& \le \widehat{C}_{\beta}\;  \calG_{e,\beta}(t+s,x),
\end{align*}
which completes the proof of Lemma \ref{L:SG-Geb}.
\end{proof}

\section{Proof of Theorem \ref{T:ExUni}}
\label{S:ExUni}
The proof of Theorem \ref{T:ExUni} will be presented at the end of this section.
Before proving Theorem \ref{T:ExUni}, we need several results.
The first one is related to the tails of the Green functions.
The corresponding results for the SHE, the SFHE, and the SWE can be found in
\cite[Proposition 5.3]{ChenDalang13Heat}, \cite[Proposition 4.7]{ChenDalang14FracHeat}, and \cite[Lemma 3.2]{ChenDalang14Wave}, respectively.
We need some notation: for $\tau>0$, $\alpha>0$ and $(t,x)\in\R_+^*\times\R$,
denote
\begin{align*}
B_{t,x,\tau,\alpha}:=\left\{(t',x')\in\R_+^*\times\R:\: 0 \le t'\le
t+\tau,\: |x-x'|\le \alpha\right\}.
\end{align*}

\begin{proposition}\label{P:G-Margin}
Suppose that $\beta\in \;]0,2[\:$.
Then for all $\tau>0$, $\alpha>0$ and $(t,x)\in\R_+^*\times\R$, there exists a constant $A>0$ such that for all
$(t',x')\in B_{t,x,1/2,1}$ and all $s \in [0,t'[$
and $y\in \R$ with
$|y|\ge A$, we have that $G_\beta\left(t'-s,x'-y\right) \le G_\beta\left(t+1-s,x-y\right)$.
\end{proposition}

\begin{proof}
Fix $(t,x) \in \R_+^*\times\R$.
By the scaling and asymptotic properties of
the Green function $G_\beta(\cdot,\cdot)$, we
know that
\begin{align*}
\frac{G_\beta(t+1-s,x-y)}{G_\beta(t'-s,x'-y)}
&=
\left(\frac{t'-s}{t+1-s}\right)^{\beta/2+1-\Ceil{\beta}}
\frac{G_\beta\left(1,\frac{x-y}{(t+1-s)^{\beta/2}}
\right)}
{G_\beta\left(1, \frac{x'-y}{(t'-s)^{\beta/2}}\right) }\\
&\approx
\left(\frac{t'-s}{t+1-s}\right)^{\frac{\beta}{2}
+a+1-\Ceil{\beta}}
\frac{|x-y|^a}{|x'-y|^a}
\;\exp\left(
\frac{b|x'-y|^c}{(t'-s)^{\beta c/2}}
-\frac{b|x-y|^c}{(t+1-s)^{\beta c/2}}
\right),
\end{align*}
as $|y|\rightarrow +\infty$ where $a=\frac{1+\beta-2\Ceil{\beta}}{2-\beta}$,
$b\in
\;]0,1[\;$ and $c>1$ (see \eqref{E:Asym-a}, \eqref{E:Asym-b} and \eqref{E:Asym-c}).
Denote
\[
f\left(\beta\right)= \frac{\beta}{2}+a +1-\Ceil{\beta}=
\frac{3+(1-\beta/2)\beta-(4-\beta)\Ceil{\beta}}{2-\beta}\;,
\]
which is plotted in Figure \ref{F6:f_beta}.
\begin{figure}[h!tbp]
 \centering
\includegraphics[scale=1]{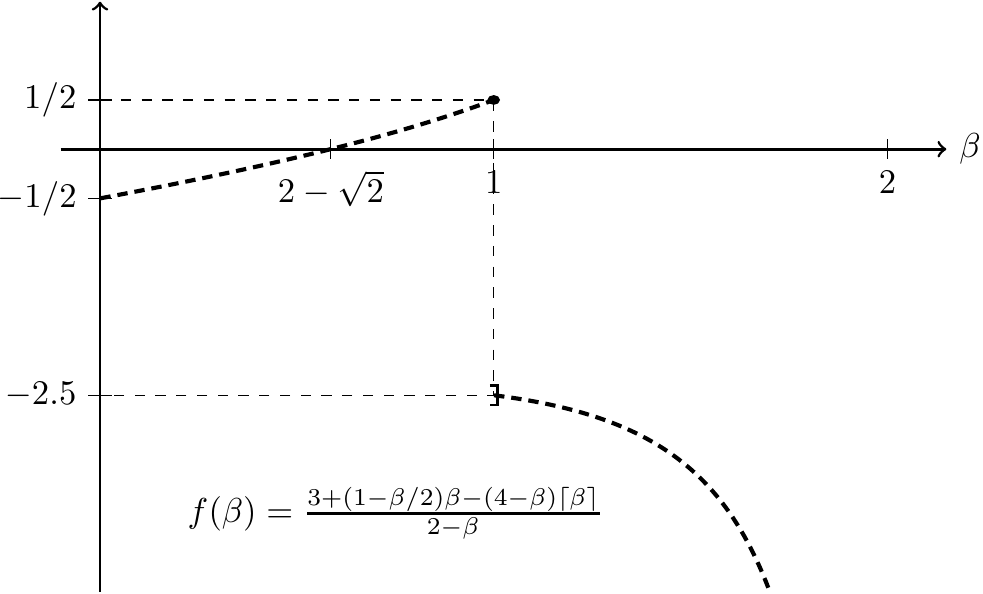}
\caption{Plot of the function $f\left(\beta\right)
$.}
\label{F6:f_beta}
\end{figure}
Let $\beta_0:=2-\sqrt{2}$.
Simple calculations show that $f\left(\beta\right)
>0$ if and only if
$\beta_0<\beta\le 1$, otherwise $f\left(\beta\right)\le 0$.
Notice that
\begin{align}\label{E_:t-s}
 \frac{t+1-s}{t'-s} = 1+\frac{t+1-t'}{t'-s} \ge
1+ \frac{t+1-t'}{t'} \ge
\frac{t+1}{t+1/2} = 1+ \frac{1}{2t +1} >1\;.
\end{align}
Hence, if $0<\beta\le \beta_0$ or $1<\beta<2$,
then
\[
\left(\frac{t'-s}{t+1-s}\right)^{f\left(\beta\right)}=
\left(\frac{t+1-s}{t'-s}\right)^{\left|f\left(\beta\right)\right|}
\ge 1\;.
\]
If $\beta_0< \beta\le 1$, we have that
\[
\left(\frac{t'-s}{t+1-s}\right)^{f\left(\beta\right)}\ge
\left(\frac{t'-s}{t+1}\right)^{\left|f\left(\beta\right)\right|}
= (t+1)^{-\left|f\left(\beta\right)\right|}
\exp\left(\left|f\left(\beta\right)\right|\;\log(t'-s)\right).
\]

Assume $n>1$. When $|y|\ge |x|+n$, we have $|x-
y|>n$ and then
\[
\frac{n}{n+1}\le \frac{|x-y|}{|x-y|+1} \le \frac{|x-y|}{|x'-y|} \le
\frac{|x-y|}{|x-y|-1}
\le \frac{n}{n-1}\;.
\]
Because $\frac{n}{n+1}>\frac{n-1}{n}$ for all
$n>1$, we have that
\begin{align}\label{E_:x-a}
\frac{|x-y|^a}{|x'-y|^a} \ge \left(1-\frac{1}{n}\right)^{|a|}\;,
\end{align}
which holds for all $a\in\R$.

The above bounds \eqref{E_:t-s} and \eqref{E_:x-a}
 imply that
\[
\left(1-\frac{1}{n}\right)\frac{|x'-y|^c}{(t'-s)^{
\beta c/2}}
\ge \left(1+\frac{1}{2t +1}\right)^{\beta c/2}
\left(1-\frac{1}{n}\right)^{c+1}
\frac{|x-y|^c}{(t+1-s)^{\beta c/2}}\:.
\]
By choosing $n$ large enough, in particular,
\[
n>\left(1-\left[1+\frac{1}{2t+1}\right]^{-\frac{\beta c}{2(c+1)}}\right)^{-1}\;,
\]
we have that
\[
\eta:= \left(1+\frac{1}{2t +1}\right)^{\beta c/2}
\left(1-\frac{1}{n}\right)^{c+1} >1\;.
\]
Thus,
\begin{align*}
\exp\left(
\frac{b|x'-y|^c}{(t'-s)^{\beta c/2}}
-\frac{b|x-y|^c}{(t+1-s)^{\beta c/2}}
\right)
&\ge \exp\left(
b(\eta-1)\frac{|x-y|^c}{(t+1-s)^{\beta c/2}}
+ \frac{b|x'-y|^c}{n(t'-s)^{\beta c/2}}
\right)\\
&\ge \exp\left(
b(\eta-1)\frac{|x-y|^c}{(t+1)^{\beta c/2}}
+ \frac{b\left(n-1\right)^c}{n(t'-s)^{\beta c/2}}
\right).
\end{align*}

Finally, if $0<\beta\le \beta_0$ or $1<\beta<2$,
then
\[
\frac{G_\beta(t+1-s,x-y)}{G_\beta(t'-s,x'-y)}
\ge \left(1-\frac{1}{n}\right)^{|a|}
\exp\left(
b(\eta-1)\frac{|x-y|^c}{(t+1)^{\beta c/2}}
+ \frac{b\left(n-1\right)^c}{n(t+1)^{\beta c/2}}
\right)\rightarrow+\infty\;,
\]
as $|y|\rightarrow+\infty$.
Hence, we can choose a large constant $A$, such
that for all $|y|\ge A$, the
inequality
\[
\frac{G_\beta(t+1-s,x-y)}{G_\beta(t'-s,x'-y)}>1
\]
holds for all $(t',x')\in B_{t,x,1/2,1}$ and $s\in [0,t']$.
If $\beta_0<\beta\le 1$, then,
\begin{align*}
\frac{G_\beta(t+1-s,x-y)}{G_\beta(t'-s,x'-y)}
\ge &
(t+1)^{-f\left(\beta\right)}
\left(1-\frac{1}{n}\right)^{|a|}
\\
&\times \exp\left(
b(\eta-1)\frac{|x-y|^c}{(t+1)^{\beta c/2}}
+ \frac{b\left(n-1\right)^c}{n(t'-s)^{\beta c/2}}
+ f\left(\beta\right) \log(t'-s)
\right).
\end{align*}
The function
\[
g(t) = \frac{C_1}{t} + C_2 \log t\;,\quad\text{for $t>0$ and $C_1, C_2>0$\;,}
\]
has its global minimum at $t_0:=C_1/C_2$:
\[
 \min_{t\in\R_+^*} g(t) = g\left(t_0\right) =C_2\left(1+\log\left(C_1/C_2\right)\right)\;,
\]
because $g'(t) = C_2 t^{-2}\left(t- t_0\right)$,
which is negative when $t<t_0$ and positive
when $t>t_0$. Hence,
\begin{align*}
 \frac{b\left(n-1\right)^c}{n(t'-s)^{\beta c/2}}
+ f\left(\beta\right) \log(t'-s)
&=
 \frac{b\left(n-1\right)^c}{n(t'-s)^{\beta c/2}}
+ \frac{2 f\left(\beta\right)}{\beta c} \log\left[
(t'-s)^{\beta c/2}\right]\\
&\ge
\frac{2 f\left(\beta\right)}{\beta c} \left(1+\log\left(\frac{2
f\left(\beta\right)
n}{b\beta c (n-1)^c}\right)\right).
\end{align*}
Therefore,
\begin{multline*}
\frac{G_\beta(t+1-s,x-y)}{G_\beta(t'-s,x'-y)}
\ge (t+1)^{-f\left(\beta\right)}
\left(1-\frac{1}{n}\right)^{|a|}
\\
\times \exp\left(
b(\eta-1)\frac{|x-y|^c}{(t+1)^{\beta c/2}}
+\frac{2 f\left(\beta\right)}{\beta c} \left(1+\log\left(\frac{2
f\left(\beta\right)
n}{b\beta c (n-1)^c}\right)\right)
\right)\rightarrow+\infty\:,
\end{multline*}
as $|y|\rightarrow +\infty$. We can choose a large
 constant $A$, such that for
all $|y|\ge A$,  all $(t',x')\in B_{t,x,1/2,1}$
and $s\in [0,t']$,
\[
\frac{G_\beta(t+1-s,x-y)}{G_\beta(t'-s,x'-y)}>1.
\]
This completes the proof of Proposition \ref{P:G-Margin}.
\end{proof}

The second set of results, Propositions \ref{P:G-SD} and \ref{P:G-FD},
give some continuity properties of the Green functions.
We need a bound of the two-parameter
Mittag-Leffler functions, which will be used in the proof of Proposition \ref{P:G-SD}.
\begin{lemma}\label{L:Bdd-MitLef}
If $0<\alpha<1$ and $\beta\ge \alpha$, then there
exists a constant $C_{\alpha,
\beta}>0$ such that
\begin{align}\label{E:Bdd-MitLef}
0< E_{\alpha,\beta}\left(-x^\alpha\right)\le
\frac{C_{\alpha,\beta}}{1+x^\alpha}
\;,\qquad\text{for all $x\ge 0$\;.}
\end{align}
\end{lemma}
\begin{proof}
Nonnegativity is due to \eqref{E:E-CM}. The upper
 bound is due to \cite[Theorem 1.6, on p. 35]{Podlubny99FDE} with $z=-x^\alpha$.
Clearly $\arg(z)=\pi$ satisfies the required condition.
\end{proof}

\begin{proposition}\label{P:G-SD}
Suppose $0<\beta<1$. Let $C_{\beta,2}$ be the universal constant in Lemma
\ref{L:Bdd-MitLef}. Then the following two properties hold:
\begin{enumerate}[(i)]
 \item For all $t>0$ and $x,y\in\R$,
\begin{align}\label{E:G-x}
\iint_{\R_+\times\R}\ud r\ud z\: \left(
G_\beta\left(t-r,x-z\right)
-G_\beta\left(t-r,y-z\right)
\right)^2=\frac{4 C_{\beta,2}}{\pi} \: t^{1-\beta} |x-y|.
\end{align}
\item For all $s,t\in\R_+^*$ with $s\le t$, and $x\in\R$,
\begin{align}\label{E:G-t1}
\int_0^s \ud r\int_\R\ud z \left(G_\beta\left(t-r,x-z\right)
-G_\beta\left(s-r,x-z\right)
\right)^2=2 C_{\beta,2}\: (t-s)^{1-\beta/2}\;,
\end{align}
and
\begin{align}\label{E:G-t2}
\int_s^t\ud r\int_\R\ud z \;
G_\beta^2\left(t-r,x-z\right)
=\frac{C_{\beta,2}}{2} (t-s)^{1-\beta/2}\;.
\end{align}
\end{enumerate}
\end{proposition}
\begin{proof}
(i) Fix $t>0$. By Plancherel's theorem and \eqref{E:G-Fourier},
the left hand side (l.h.s.) of \eqref{E:G-x} is equal to
\begin{align*}
\frac{1}{2\pi}\int_0^t\ud r\int_\R\ud \xi &\left|
e^{-i \xi x}E_{\beta,1}\left(-(t-r)^\beta \xi^2\right)
-
e^{-i \xi y}E_{\beta,1}\left(-(t-r)^\beta
\xi^2\right)
\right|^2\\
&=\frac{1}{2\pi}\int_0^t \ud r\int_\R \ud \xi \:
E_{\beta,1}^2
\left(-(t-r)^\beta \xi^2\right)\;
\left|e^{-i\xi x}-e^{-i\xi y}\right|^2\\
&=\frac{1}{\pi} \int_\R \ud \xi\:
\left(1-\cos(\xi(x-y))\right) \int_0^t\ud r
\;E_{\beta,1}^2\left(-(t-r)^\beta \xi^2\right).
\end{align*}
By \eqref{E:E-CM},
\begin{align*}
\int_0^t\ud r
\;E_{\beta,1}
^2\left(-(t-r)^\beta \xi^2\right)
&\le
E_{\beta,1}(0) \int_0^t\ud r
\;E_{\beta,1}
\left(-(t-r)^\beta \xi^2\right)=
t\: E_{\beta,2}\left(-t^\beta \xi^2\right),
\end{align*}
where the last equality can be obtained by integration term-by-term (see also
\cite[(1.99), on p. 24]{Podlubny99FDE}).
Then use the bound \eqref{E:Bdd-MitLef} and the fact that $1-\cos(x)\le 2\wedge \left(x^2/2\right)$ for all $x\in\R$
to see that the l.h.s. of \eqref{E:G-x} is bounded by
\begin{align*}
\frac{C_{\beta,2}t^{1-\beta}}{\pi} \int_\R\ud \xi\:\frac{2\wedge
\left[(x-y)\xi\right/\sqrt{2}\:]^2}{\xi^2}
&=\frac{\sqrt{2} \: C_{\beta,2} t^{1-\beta}}{\pi}
|x-y|
\int_0^\infty\ud u\: \frac{2\wedge u^2}{u^{2}} =\frac{4 C_{\beta,2} t^{1-\beta}}{\pi}
|x-y|.
\end{align*}

(ii) Denote the l.h.s. of \eqref{E:G-t1} by $I$.  Apply Plancherel's theorem and use \eqref{E:E-CM},
\begin{align*}
 I=&\frac{1}{2\pi}\int_0^s\ud r\int_\R \ud \xi \left|
e^{-i\xi x} E_{\beta,1}\left(-(t-r)^\beta \xi^2
\right)
-
e^{-i\xi x}
 E_{\beta,1}\left(-(s-r)^\beta \xi^2\right)
\right|^2 \\
=&
\frac{1}{2\pi}\int_0^s\ud r\int_\R \ud \xi \left|
E_{\beta,1}\left(-(t-r)^\beta \xi^2\right)
-
 E_{\beta,1}\left(-(s-r)^\beta \xi^2\right)
\right|^2 \\
 \le&
\frac{1}{2\pi}\int_0^s\ud r\int_\R \ud \xi\;
2 E_{\beta,1}(0)
\left[
E_{\beta,1}\left(-(t-r)^\beta \xi^2\right)
-
 E_{\beta,1}\left(-(s-r)^\beta \xi^2\right)
\right].
\end{align*}
Integration term-by-term gives that
\begin{align*}
\int_0^s \ud r\:  E_{\beta,1}\left(-(t-r)^\beta\xi^2\right)
&=\int_{t-s}^t \ud r\: E_{\beta,1}\left(-r^\beta\xi^2\right)\\
&=t \:E_{\beta,2}\left(-t^\beta \xi^2\right)
-(t-s)\: E_{\beta,2}\left(-(t-s)^\beta \xi^2\right
)\;.
\end{align*}
Hence, by \eqref{E:E-CM} again,
\begin{align*}
I & \le
\frac{1}{\pi}
\int_\R \ud \xi \left(
t E_{\beta,2}\left(-t^\beta \xi^2 \right)
-s E_{\beta,2}\left(-s^\beta \xi^2 \right)
-(t-s) E_{\beta,2}\left(-(t-s)^\beta \xi^2 \right)
\right)\\
&\le
\frac{1}{\pi}(t-s) \int_\R\ud \xi \left(
E_{\beta,2}\left(-t^\beta\xi^2\right)+
E_{\beta,2}\left(-(t-s)^\beta\xi^2\right)
\right).
\end{align*}
Then by the bound in \eqref{E:Bdd-MitLef} and the
integral $\int_\R \ud
\xi\:\frac{1}{1+c^2 \xi^2}  =\pi/|c|$ for $c\ne 0$,
we find that
\[
I \le C_{\beta,2}(t-s)
\left(\frac{1}{t^{\beta/2}}+\frac{1}{(t-s)^{\beta/2}}\right)
\le
2 C_{\beta,2} (t-s)^{1-\beta/2}.
\]
As for \eqref{E:G-t2}, by a similar reasoning, we have
\begin{align*}
\int_s^t\ud r\int_\R\ud z \;&
G_\beta^2\left(t-r,x-z\right)
\le
\frac{1}{2\pi}\int_s^t\ud r\int_\R\ud \xi \;
E_{\beta,1}^2\left(-(t-r)^\beta\xi^2\right)\\
&\le
\frac{1}{2\pi}\int_s^t\ud r\int_\R\ud \xi \;
E_{\beta,1}\left(-(t-r)^\beta\xi^2\right)
\le
\frac{t-s}{2\pi} \int_\R\ud \xi \;
 E_{\beta,2}\left(-(t-s)^\beta\xi^2\right)\\
&\le
\frac{t-s}{2\pi} \int_\R \ud \xi\;
 \frac{C_{\beta,2}}{1+(t-s)^\beta\xi^2}
=\frac{C_{\beta,2}}{2} (t-s)^{1-\beta/2}\;,
\end{align*}
which completes the proof of Proposition \ref{P:G-SD}.
\end{proof}

For the fast diffusion equations, we are only able to prove the following less precise results
in Proposition \ref{P:G-FD}  due to the lack of complete monotonicity for
$E_{\alpha,\beta}(-x)$ with $\alpha>1$; see \eqref{E:E-CM} for the necessary
and sufficient conditions for $E_{\alpha,\beta}(-x)$ to be completely monotonic.

\begin{proposition}\label{P:G-FD}
For all $(t,x)\in\R_+\times\R$ and $1<\beta<2$,
we have
\[
\lim_{(t',x')\rightarrow (t,x)} \iint_{\R_+\times\R}\ud s\ud y\: \left(
G_\beta\left(t'-s,x'-y\right)
-G_\beta\left(t-s,x-y\right)
\right)^2=0.
\]
\end{proposition}
\begin{proof}
We only need to consider the case where $t>0$.
Fix $(t,x)\in\R_+^*\times\R$.
Denote $\Lambda:=\sup_{x\in\R} G_\beta(1,x)$.
We are going to apply the Lebesgue dominated convergence theorem. Clearly, by
the continuity of the Green functions, for all $(s,y)\in\R_+^*\times\R$,
\[
G_\beta\left(t'-s,x'-y\right)
-G_\beta\left(t-s,x-y\right) \rightarrow 0\;,\quad \text{as
$(t',x')\rightarrow (t,x)$.}
\]
 We need to find
an integrable bound.
Choose $A>0$ according to Proposition \ref{P:G-Margin} and suppose that
$(t',x')\in B_{t,x,1/2,1}$.
If $|y|>A$, since $1-\beta/2>0$, by
Proposition \ref{P:G-Margin},
\begin{align*}
\left|
G_\beta\left(t'-s,x'-y\right)
-G_\beta\left(t-s,x-y\right)
\right|^2
&\le
4 G_\beta^2\left(t+1-s,x-y\right)\\
&\le 4\Lambda(t+1-s)^{1-\beta/2}
G_\beta\left(t+1-s,x-y\right)\\
&\le 4\Lambda(t+1)^{1-\beta/2}
G_\beta\left(t+1-s,x-y\right).
\end{align*}
If $|y|\le A$, we have that
\begin{align*}
\left|
G_\beta\left(t'-s,x'-y\right)
-G_\beta\left(t-s,x-y\right)
\right|^2
&\le
2 G_\beta^2\left(t'-s,x'-y\right)+2 G_\beta^2\left(t-s,x-y\right)\\
&\le 2\Lambda^2 \left[(t'-s)^{2-\beta}+(t-s)^{2-\beta}\right]\\
&\le 4\Lambda^2 (t+1)^{2-\beta}\;.
\end{align*}
Hence,
\begin{multline*}
\left|
G_\beta\left(t'-s,x'-y\right)
-G_\beta\left(t-s,x-y\right)
\right|^2\\
\le 4\Lambda(t+1)^{1-\beta/2}
G_\beta\left(t+1-s,x-y\right)\Indt{|y|>A} +
4\Lambda^2 (t+1)^{2-\beta} \Indt{|y|<A,\;0\le s\le t+1}\;.
\end{multline*}
Denote this upper bound by $f(s,y)$.
Clearly, this upper bound is integrable:
\begin{align*}
\int_{\R_+}\ud s\int_\R \ud y\: f(s,y)
\le&
4\Lambda(t+1)^{1-\beta/2} \iint_{\R_+\times\R} \ud s\ud y\:
G_\beta\left(t+1-s,x-y\right)\\
&+
4\Lambda^2 (t+1)^{2-\beta}
\iint_{[0,t+1]\times[-A,A]}\ud s\ud y\\
=&
2\Lambda (t+1)^{3-\beta/2} +
8 A \Lambda^2 (t+1)^{3-\beta}<+\infty\:.
\end{align*}
Therefore, this proposition is proved by the
Lebesgue dominated convergence theorem.
\end{proof}

The third result, Proposition \ref{P:J0LocalLip}, is about solutions to the homogeneous equation.
We need to prove a lemma first.
Recall the function $f_\eta(x)$ defined in  \eqref{E:f}.

\begin{lemma}\label{L:J0Lip}
Suppose $\beta \in\;]0,2[\;$.
Let $b\in \:]0,1[\;$ be the constant defined in \eqref{E:Asym-b}.
Then for all $\eta>0$, the following three functions
\[
f_\eta\left(\frac{x}{t^{\beta/2}}\right)\;,\quad
G_\beta(t,x) \;f_b^{-1}\left(\frac{x}{t^{\beta/2}}
\right)\;,
\quad\text{and}\quad
G_\beta^*(t,x) \;f_b^{-1}\left(\frac{x}{t^{\beta/2}}\right)
\]
are Lipschitz continuous over
$\R_+^*\times\R$, that is, for all $x,y\in\R$ and
$t,s\ge \epsilon >0$, there
exists a constant $C_\epsilon>0$ such that
\begin{gather}\label{E:f-Lip}
\left|
f_\eta\left(\frac{x}{t^{\beta/2}}\right)
-
f_\eta\left(\frac{y}{s^{\beta/2}}\right)
\right|
\le C_\epsilon \left(|x-y|+|t-s|\right)\;,\\
\label{E:Gf-Lip}
\left|
G_\beta(t,x)f_b^{-1}\left(\frac{x}{t^{\beta/2}}\right)
-
G_\beta(s,y)f_b^{-1}\left(\frac{y}{s^{\beta/2}}\right)
\right|
\le C_\epsilon \left(|x-y|+|t-s|\right)\;,\\
\label{E:G*f-Lip}
\left|
G_\beta^*(t,x)f_b^{-1}\left(\frac{x}{t^{\beta/2}}\right)
-
G_\beta^*(s,y)f_b^{-1}\left(\frac{y}{s^{\beta/2}}\right)
\right|
\le C_\epsilon \left(|x-y|+|t-s|\right).
\end{gather}
\end{lemma}
\begin{proof}
(i) We first prove \eqref{E:f-Lip}.
Denote
\[
g(t,x) = f_\eta\left(\frac{x}{t^{\beta/2}}\right)
=
\exp\left(-\frac{\eta}{2 t^{\beta c/2}}|x|^c\right
)\;,
\]
where $c=\frac{2}{2-\beta}$. Fix $x\ne 0$. Clearly
,
\[
\left|\frac{\partial }{\partial x}g(t,x)\right| =
\frac{\eta \; c |x|^{c-1}}{2\; t^{\beta c/2}}
g(t,x)\le \frac{c\; \eta}{2\; t^{\beta/2}} \sup_{y\in\R} |y|^{c-1}f_\eta(y)\;,
\]
and
\[
\left|\frac{\partial }{\partial t}g(t,x)\right| =
\frac{\eta \beta c |x|^c}{4\; t^{\beta c/2-1}}
g(t,x)
\le
\frac{\eta\;\beta \;c}{4 \;t} \sup_{y\in\R} |y|^c
f_\eta(y)\;.
\]
Note that the two suprema are finite because $c>1$.
Hence, by the mean value theorem,
\begin{align*}
\left|g(t,x)-g(s,y)\right|
&\le \left|g(t,x)-g(t,y)\right| + \left|g(t,y)-g(
s,y)\right|\\
&\le \frac{C_1}{\epsilon^{\beta/2}} \left|x-y\right|
+\frac{C_2}{\epsilon} |t-s|
\end{align*}
for $xy\ge 0$ (i.e., $x$ and $y$ have the same
sign) and $t,s\ge \epsilon>0$,
where
\[
C_1 = \frac{\eta}{2-\beta}
\sup_{y\in\R}|y|^{\frac{\beta}{2-\beta}}f_\eta(y)\quad\text{and}\quad
C_2=\frac{\eta \beta}{2(2-\beta)} \sup_{y\in\R}
|y|^{\frac{2}{2-\beta}}
f_\eta(y)\;.
\]
When $x$ and $y$ have different signs, we use the
fact that
\begin{align}\notag
\left|g(t,x)-g(t,y)\right|
&\le
\left|g(t,x)-g(t,0)\right|+
\left|g(t,0)-g(t,y)\right|\\
&\le
\frac{C_1}{\epsilon^{\beta/2}}
\left(|x|+|y|\right)
=
\frac{C_1}{\epsilon^{\beta/2}}
\left|x-y\right|.
\label{E:xy<0}
\end{align}

(ii) Now let us prove \eqref{E:Gf-Lip}.
Denote
\[
h(t,x) = G_\beta(t,x)f_b^{-1}\left(\frac{x}{t^{\beta/2}}\right).
\]
We assume that both $x$ and $y$ are nonnegative.
The case where $xy<0$ can be covered by a similar argument as \eqref{E:xy<0}.
By \eqref{E:G-dx},
\begin{align*}
\frac{\partial }{\partial x} h(t,x) =&
-\frac{t^{\Ceil{\beta}-1-\beta}}{2}M_{\beta/2,\Ceil{\beta}-\beta/2}
\left(\frac{x}{t^{\beta/2}}\right)
\exp\left(\frac { b } { 2 } \left|\frac{x}{ t^ {\beta/2 } }\right|^c\right)
\\
&+
\frac{b c t^{\Ceil{\beta}-1-\beta}}{4}
\left|\frac{x}{t^{\beta/2}}\right|^{c-1}
\exp\left(\frac{b}{2}\left|\frac{x}{t^{\beta/2}}\right|^c\right)
M_{\beta/2,\Ceil{\beta}}
\left(\frac{x}{t^{\beta/2}}\right).
\end{align*}
Because the exponent of the asymptotic of $M_{\lambda,\theta}(x)$ in
\eqref{E:Asy-M} depends only on the first parameter $\lambda$, we have that
\[
C_{\alpha,\theta}:= \sup_{y\in\R} |y|^{\alpha} \left|
 M_{\beta/2, \theta}(y)\right|\:
\exp\left(\frac{b}{2} |y|^c\right)< +\infty\;,\quad\text{for all $\theta\in\R$ and
$\alpha>0$\:.}
\]
Note that $\Ceil{\beta}-1-\beta\le 0$.
Therefore,
\[
\left|\frac{\partial}{\partial x}h(t,x)\right|
\le C_\epsilon'\;,\quad\text{for all $x\in\R$ and
$t\ge \epsilon>0$,}
\]
where
\[
C_\epsilon' =
\frac{\epsilon^{\Ceil{\beta}-1-\beta}}{2} C_{0,\Ceil{\beta}-\beta/2} +
\frac{b\; c\; \epsilon^{\Ceil{\beta}-1-\beta}}{4}
C_{c-1,\Ceil{\beta}}\;.
\]
By \eqref{E:Der-M}, we have
\begin{align*}
\frac{\partial}{\partial t}G_\beta(t,x)
=&\frac{\Ceil{\beta}-1-\beta/2}{2} t^{\Ceil{\beta}
-2-\beta/2}
M_{\beta/2,\Ceil{\beta}}\left(\frac{x}{t^{\beta/2}}\right)
\\
&+\frac{\beta}{2} \left(\frac{x}{t^{\beta/2}}\right)
t^{\Ceil{\beta}-2-\beta/2}M_{\beta/2,\Ceil{\beta}-\beta/2}
\left(\frac{x}{t^{\beta/2}}\right).
\end{align*}
Notice that
\[
\frac{\partial}{\partial t}
\exp\left(\frac { b } { 2 } \left|\frac{x}{ t^{\beta/2 } }\right|^c\right)
=
-\frac{\beta \; c \; b}{4} \left(\frac{x}{t^{\beta/
2}}\right)^c t^{-1}
\exp\left(\frac { b } { 2 } \left|\frac{x}{ t^ { \beta/2 } }\right|^c\right).
\]
Hence,
\begin{align*}
 \frac{\partial }{\partial t} h(t,x)
=&\frac{\Ceil{\beta}-1-\beta/2}{2} t^{\Ceil{\beta}-
2-\beta/2}
M_{\beta/2,\Ceil{\beta}}\left(\frac{x}{t^{\beta/2}
}\right)
\exp\left(\frac { b }{ 2 } \left|\frac{x}{ t^{\beta/2 } }\right|^c\right)
\\
&+\frac{\beta}{2} \frac{x}{t^{\beta/2}}
t^{\Ceil{\beta}-2-\beta/2}M_{\beta/2,\Ceil{\beta}-
\beta/2}
\left(\frac{x}{t^{\beta/2}}\right)\exp\left(\frac{
 b } { 2 } \left|\frac{x}{
t^{\beta/2 } }\right|^c\right) \\
&-\frac{\beta \; c \; b}{8} \left(\frac{x}{t^{\beta/2}}\right)^c
t^{\Ceil{\beta}-2-\beta/2}
M_{\beta/2,\Ceil{\beta}}\left(\frac{x}{t^{\beta/2}
}\right)\exp\left(\frac { b }
{ 2 } \left|\frac{x}{ t^ { \beta/2 } }\right|^c\right).
\end{align*}
Note that $\Ceil{\beta}-2-\beta/2\le 0$. Therefore,
\[
\left|\frac{\partial}{\partial x}h(t,x)\right|
\le C_\epsilon''\;,\quad\text{for all $x\in\R$
and $t\ge \epsilon>0$,}
\]
where
\[
C_\epsilon'' =\epsilon^{\Ceil{\beta}-2-\beta/2}\left(
\frac{\left|\Ceil{\beta}-1-\beta/2\right|}{2} C_{0,\Ceil{\beta}} +
\frac{\beta}{2}
C_{1,\Ceil{\beta}-\beta/2}
+\frac{\beta\; c\; b}{8}
C_{c,\Ceil{\beta}}
\right).
\]
Finally, apply the mean value theorem to conclude this case.
The argument is the same as (i).

(iii) \eqref{E:G*f-Lip} can be proved in the same way.
We will not repeat here.
This completes the proof of Lemma \ref{L:J0Lip}.
\end{proof}

\begin{proposition}\label{P:J0LocalLip}
Suppose that $0<\beta<2$ and $\mu\in\calM_T^{\beta}\left(\R\right)$. Denote
\[
J_1(t,x):=\left(G_\beta(t,\cdot) *\mu\right)(x)\quad \text{and}\quad
J_2(t,x):=\left(G_\beta^*(t,\cdot) *\mu\right)(x).
\]
(1) Both functions $J_i(t,x)$, $i=1,2$, are
locally Lipschitz continuous on $\R_+^*\times\R$,
that is,
for all
compact sets $K\subseteq \R_+^*\times\R$, there
exits a constant $C_K>0$
such that
\[
\left|J_i(t,x) - J_i(s,y)\right| =
C_K \left(|t-s|+|x-y|\right)\;,\quad\text{for all $(t,x)$ and $(s,y)\in K$.}
\]
Hence, the solution $J_0(t,x)$ in \eqref{E:J0} is locally Lipschitz continuous on $\R_+^*\times\R$.\\
(2) If $0<\beta\le 1$ and if $\mu(\ud x) =f(x)\ud x$ where $f$ is
$\alpha$--H\"older continuous with $\alpha\in\;]0,
1]$, then $J_1(\cdot,\circ) \in
C_{\alpha\beta/2,\:\alpha\beta}\left(\R_+\times\R\right)$.
\end{proposition}
A similar proof for part (2) for SHE can be found in \cite[Lemma 3.8]{ChenDalang13Heat}.
\begin{proof}
(1) We first show the Lipschitz continuity of
the function
 $(t,x)\mapsto J_1(t,x)$ for $0<\beta<2$.
Let $\epsilon = \inf\left\{s: (s,y)\in K\right\}$,
 $T = \sup\left\{s:
(s,y)\in K\right\}$ and $k= \sup\left\{|y|: (s,y)\in K \right\}$.
Since $K$ is a compact set of $\R_+^*\times\R$,
we know that $\epsilon>0$,
$T<+\infty$ and $k<+\infty$.
Suppose $\epsilon \le t,s\le T$ and $x,y \in [-k,k]$.
Notice that
\[
\left|J_1(t,x) - J_1(s,y)\right|  \le
\int_\R  |\mu|(\ud z) \:\left|G_\beta(t,x-z) -G_\beta(s,y-z) \right|
\]
and
\begin{align*}
\big|G_\beta(t,x-z) -&G_\beta(s,y-z) \big|\\
\le
&\quad \left|G_\beta(t,x-z)f_b^{-1}\left(\frac{x-z}{t^{\beta/2}}\right) -G_\beta(s,y-z)
f_b^{-1}\left(\frac{y-z}{s^{\beta/2}}\right)
\right|f_b\left(\frac{x-z}{t^{\beta/2}}\right)\\
&+
\left|
f_b\left(\frac{x-z}{t^{\beta/2}}\right)-
f_b\left(\frac{y-z}{s^{\beta/2}}\right)
\right| G_\beta(s,y-z) f_b^{-1}\left(\frac{y-z}{s^{\beta/2}}\right),
\end{align*}
where $b>0$ is defined in \eqref{E:Asym-b}.
By Lemma \ref{L:J0Lip}, there is a constant $C_\epsilon>0$ such that
\[
\left|G_\beta(t,x-z) -G_\beta(s,y-z) \right| \le
C_\epsilon \left(|t-s|+|x-y|\right).
\]
By the asymptotics of $G_\beta(s,y)$ with $s$
fixed, we know that for some
constant $C>0$,
\begin{align}\label{E_:J0LocalLip}
G_\beta(s,y-z) f_b^{-1}\left(\frac{y-z}{s^{\beta/2}}\right)
\le
C \; s^{\Ceil{\beta}-1-\beta/2} f_{b/2}\left(\frac{y-z}{s^{\beta/2}}\right)
\;,\quad\text{for all $z\in\R$\;.}
\end{align}
Notice that $f_b\left(\frac{x}{t^{\beta/2}}\right)
 \le
f_{b \;\epsilon^{\beta c/2}}(x)$ if $t\ge \epsilon
$,
where $c=\frac{2}{2-\beta}$. Since
\[
-\frac{1}{2}\le \Ceil{\beta}-1-\beta/2\le \frac{1}{2}\;,
\]
we have
\[
s^{\Ceil{\beta}-1-\beta/2} \le \sqrt{T \vee \epsilon^{-1}}\:,\quad \text{for $\epsilon\le s\le T$}\:,
\]
where $a\vee b:=\max(a,b)$.
Therefore,
\[
\left|J_1(t,x) - J_1(s,y)\right|  \le
C_\epsilon
\left[
\left(|\mu|*f_{b \;\epsilon^{\beta c/2}}\right)(x)
+
C\sqrt{T \vee \epsilon^{-1}}\left(|\mu|*f_{b \;\epsilon^{\beta
c/2}/2}\right)(y)
\right]
\left(|t-s|+|x-y|\right),
\]
for all $x,y\in\R$ and $t,s\ge \epsilon$.
The function $x\mapsto \left(|\mu|*f_\eta\right)(x)$ is well defined because
$\mu\in\calM_T^\beta(\R)$. Moreover, it is continuous, which can be easily proved by the dominated convergence theorem thanks  to the continuity and
boundedness of $f_\eta(x)$.

As for the function $J_2(t,x)$, we simply change
the power of $s$ in
\eqref{E_:J0LocalLip} by $-\beta/2$ and so
\[
s^{-\beta/2} \le \epsilon^{-\beta/2}\le \epsilon^{-1}\;,\quad\text{for $s\ge
\epsilon$ and $1<\beta<2$}\;.
\]
Hence, we need to replace the term $\sqrt{T\vee\epsilon^{-1}}$ by $\epsilon^{-1}$.
Clearly, $\sqrt{T\vee \epsilon^{-1}} \vee \epsilon^{-1} = \sqrt{T}\vee \epsilon^{-1}$.
Finally, we can choose the following constant for both $J_1(t,x)$ and $J_2(t,x)$:
\[
C_K=2 C_\epsilon
\left(
\sup_{x\in [-k,k]}\left(|\mu|*f_{b \;\epsilon^{\beta c/2}}\right)(x)+
C\left(\sqrt{T}\vee \epsilon^{-1}\right)\sup_{x\in
[-k,k]}\left(|\mu|*f_{b \;\epsilon^{\beta c/2}/2}\right)(x)
\right)<+\infty\;.
\]

(2) Fix $(t,x)$ and $(t',x')\in\R_+\times\R$ with $t'>t$. Then we have that
\begin{align*}
\left|J_1(t,x)-J_1(t',x')\right|
&\le
\left|J_1(t,x)-J_1(t',x)\right|+\left|J_1(t',x)-J_1(t',x')\right|\\
&:=I_1(t,t';x)+I_1(t';x,x').
\end{align*}
By change of variables and the H\"older continuity of $f$, for some constant $C>0$,
\begin{align*}
 I_1(t,t';x) &= \left|\int_\R\ud y
\:\left(G_\beta(t,x-y)-G_\beta(t',x-y)\right) f(y)
\right|\\
&=
\left|\int_\R\ud z
\:G_\beta(1,z) \left(f(x-t^{\beta/2}\:z)-f(x-(t')
^{\beta/2}\:z)\right)
\right|\\
&\le C\left|t^{\beta/2}-(t')^{\beta/2}\right|^\alpha \int_\R\ud z\:
G_\beta(1,z)|z|^\alpha,
\end{align*}
where the integral is finite by \eqref{E:G-Mom}.
By subadditivity of the function $x\in\R_+\mapsto x^{\beta/2}$,
$(t')^{\beta/2}-t^{\beta/2}\le |t'-t|^{\beta/2}$.
The arguments for $I_2(t';x,x')$ are similar. We will not repeat here.
This completes the proof of Lemma \ref{P:J0LocalLip}.
\end{proof}

The last result, Lemma \ref{L:InDt}, is about the initial data.
Similar results for the SHE, the SFHE, and the SWE can be proved in
\cite[Lemma 3.9]{ChenDalang13Heat}, \cite[Lemma 4.9]{ChenDalang14FracHeat} and \cite[Lemma 3.4]{ChenDalang14Wave}, respectively.
Recall that $J_0(t,x)$ is the solution to the
homogeneous equation; see \eqref{E:J0}.

\begin{lemma}\label{L:InDt}
Suppose $0<\beta<2$. For all $\mu$ and $\nu\in \calM_T^\beta\left(\R\right)$, all
compact sets $K\subseteq \R_+^*\times\R$,
\[
\sup_{(t,x)\in K}\left(\left[1+J_0^2\right]\star \calK\right)(t,x)
<\infty.
\]
\end{lemma}
\begin{proof}
We need only consider the part $J_0^2\star \calK$ because the part
$1\star \calK$ can be obtained by the special case where $\mu(\ud x) = \ud x$.
Assume that $\mu \ge 0$. For general $\mu$, we simply replace $\mu$ below
by $|\mu|$.
The case $\beta=1$ is covered by
\cite[Lemma 3.9]{ChenDalang13Heat}.
Note that by \eqref{E:UpBd-K-A} and Proposition \ref{P:ST-Con}, for two
constants $c_1$ and $c_2>0$, one has that for all
$(t,x)\in\R_+\times\R$,
\begin{align}
 \label{E_:h}
\calK(t,x)\le \G_\beta(t,x) h(t),\quad \text{with}\quad
h(t):= c_1 \left(t^{-\sigma}+
e^{c_2\:t}\right),
\end{align}
where $\G_\beta(t,x)$ is defined in \eqref{E:calG} and $\sigma<1$ is defined in \eqref{E:SigUps}.
In the following, denote $\bar{z}=(z_1+z_2)/2$ and $\Delta z=z_1-z_2$.

\paragraph{Slow diffusions}  Fix $\beta\in \:]0,1[\:$.
By the same argument as Proposition \ref{P:ST-Con}, for some nonnegative constant
$C_\beta<\infty$,
$G_\beta(t,x)\le C_\beta \: \calG_{e,\beta}(t,x)$ for all $(t,x)
 \in\R_+\times\R$.
Thus,
\[
J_0^2(s,y) \le C_\beta^2 \iint_{\R^2} \mu(\ud z_1)
\mu(\ud z_2)
\: \calG_{e,\beta}(s,y-z_1)\calG_{e,\beta}(s,y-z_2).
\]
Because $K(t,x)\le C_\beta \: \calG_{e,\beta}(t,x)\: h(t)$, we
see that
\begin{align*}
\left(J_0^2 \star \calK\right)(t,x)
\le& C_\beta^2 \int_0^t \ud s \: h(t-s) \int_\R \ud
y \iint_{\R^2}
\mu(\ud z_1)\mu(\ud z_2)
\\
&\times \calG_{e,\beta}(s,y-z_1)\; \calG_{e,\beta}(s,y-z_2)\; \calG_{e,\beta}(t-s,x-y).
\end{align*}
By the inequality
\begin{align}\label{E:BtDz}
\left|\frac{a+b}{2}\right|+\left|\frac{a-b}{2}\right| \le |a|+ |b|,
\end{align}
we see that $\calG_{e,\beta}(s,y-z_1)\; \calG_{e,\beta}(s,y-z_2)\le \calG_{e,\beta}\left(s,
y-\bar{z}\right)
\;\calG_{e,\beta}\left(s, 1/2\:\Delta z\right)$.
Then integrate over $\ud y$ using Lemma \ref{L:SG-Geb},
\begin{align*}
\left(J_0^2 \star \calK\right)(t,x)
\le C_\beta^2 \: \widehat{C}_\beta \int_0^t \ud s \: h(t-s) \iint_{\R^2}
\mu(\ud z_1)\mu(\ud z_2)\: \calG_{e,\beta}(s,1/2 \:\Delta z)\;
\calG_{e,\beta}\left(2^{1/\beta}t,x-\bar{z}\right).
\end{align*}
By \eqref{E:BtDz} again,
\begin{align*}
\calG_{e,\beta}(s,1/2\:\Delta z)\;
\calG_{e,\beta}\left(2^{1/\beta}t,x-\bar{z}\right)&=
\frac{1}{4\sqrt{2} \:(s t)^{\beta/2}}
\exp\left(-\frac{|z_1-z_2|}{2 \: s^{\beta/2}}
-\frac{\left|(x-z_1)+(x-z_2)\right|}{\sqrt{8} \:
t^{\beta/2}}
\right)\\
&\le
\frac{1}{4\sqrt{2} \:(s t)^{\beta/2}}
\exp\left(-\frac{|z_1-z_2|}{\sqrt{8} \: t^{\beta/2}}
-\frac{\left|(x-z_1)+(x-z_2)\right|}{\sqrt{8} \:
t^{\beta/2}}
\right)\\
&
\le \frac{1}{4\sqrt{2} \:(s t)^{\beta/2}}
\exp\left(-\frac{|x-z_1|}{\sqrt{8} \: t^{\beta/2}}
-\frac{\left|x-z_2\right|}{\sqrt{8} \: t^{\beta/2}
}
\right)\\
&= 4\sqrt{2} \: \calG_{e,\beta}\left(2^{6/\beta} t, x-z_1\right)
\:
\calG_{e,\beta}\left(2^{6/\beta} t, x-z_2\right).
\end{align*}
Denote $I(t,x)= \left(\mu * \calG_{e,\beta}\left(2^{6/\beta} t,\cdot\right)\right)(x)$.
Clearly, $\mu\in
\calM_T^\beta(\R)$ implies that $I(t,x)<+\infty$
for all
$(t,x)\in\R_+^*\times\R$.
Therefore,
\[
\left(J_0^2 \star \calK\right)(t,x)\le
4\sqrt{2}\: C_\beta^2 \: \widehat{C}_\beta\:  I(t,x)^2  \int_0^t \ud s
\: h(t-s).
\]
Clearly, the $\ud s$-integral is integrable because $\sigma<1$.

\paragraph{Fast diffusions} Fix $\beta\in \:]1,2[\:$.
By \eqref{E:J0}, we only need  to consider two  cases: Case I --- $\mu= 0$ and $\nu\ne 0$, and Case II --- $\mu\ne 0$ and $\nu= 0$.

We first consider Case I.
By the same arguments as in Proposition \ref{P:ST-Con}, for some nonnegative constant $C_\beta<+ \infty$,
$G_\beta(t,x)\le C_\beta t\: G_1(t^\beta,x)$. Thus,
\[
J_0^2(s,y)\le C_\beta^2 s^2 \iint_{\R^2}\nu(\ud z_1)\nu(\ud
z_2)\:
G_1(s^\beta,y-z_1)G_1(s^\beta,y-z_2)
\]
By \eqref{E_:h},  we see that
\begin{align*}
\left(J_0^2 \star \calK\right)(t,x)
\le &C_\beta^2 \int_0^t \ud s \: h(t-s) s^2 \int_\R \ud y
\iint_{\R^2}
\nu(\ud z_1)\nu(\ud z_2)
\\
& \times G_1(s^\beta,y-z_1)\; G_1(s^\beta,y-z_2)\;
G_1((t-s)^\beta,x-y).
\end{align*}
 By \cite[Lemma 5.4]{ChenDalang13Heat},
\[
G_1(s^\beta,y-z_1)\; G_1(s^\beta,y-z_2)
= G_1\left(\frac{s^\beta}{2},y- \bar{z}\right)  \:
G_1\left(2 s^\beta, \Delta z \right)
\le \sqrt{2}\:  G_1\left(s^\beta,y- \bar{z}\right)
G_1\left(2 s^\beta, \Delta z \right).
\]
Integrate over $\ud y$ using the semigroup property of the heat kernel
function,
\[
\left(J_0^2 \star \calK\right)(t,x)
\le
C_\beta^2\sqrt{2} \int_0^t \ud s \: h(t-s) s^2\iint_
{\R^2}
\nu(\ud z_1)\nu(\ud z_2)
G_1(2s^\beta,\Delta z) G_1((t-s)^\beta+s^\beta,x-\bar{z}).
\]
By \eqref{E:t-s+s_F} and \cite[Lemma 5.5]{ChenDalang13Heat},
\begin{align*}
G_1((t-s)^\beta+s^\beta,x-\bar{z})\: G_1(2s^\beta,
\Delta z)
& \le
2^{\frac{\beta-1}{2}} G_1(t^\beta, x-\bar{z})\: G_1(2s^\beta,\Delta z)\\
& \le 2^{1+\beta/2} \frac{t^{\beta/2}}{s^{\beta/2}
}
G_1(4t^\beta,x-z_1)G_1(4t^\beta,x-z_2).
\end{align*}
Denote $I(t,x) = \int_\R \nu(\ud z)\: G_1(4t^\beta,x-z)$. Clearly, $\nu\in\calM_T^\beta(\R)$ implies that $I(t,x)<+\infty$
for all
$(t,x)\in\R_+^*\times\R$.
Therefore,
\[
\left(J_0^2 \star \calK\right)(t,x)\le
C_\beta^2 \: 2^{\frac{\beta +3}{2}} t^{\beta/2} I(t,x)^2  \int_0^t \ud s \:s^{2-\beta/2} h(t-s).
\]
Clearly, the $\ud s$-integral is integrable because $\sigma<1$ and $2-\beta/2>-1$.

As for Case II, by the same argument as Proposition \ref{P:ST-Con}, for some nonnegative
constant
$C_\beta<\infty$,
$G_\beta^*(t,x)\le C_\beta \: g(t,x)$ for all $(t,x) \in\R_+\times\R$.
Therefore, this case can be proved by the same
arguments as the slow diffusion case with $G_\beta(t,x)$ replaced by $G_\beta^*(t,x)$.

Finally, we remark that in both cases, by the continuity of the function
$\R_+^*\times\R \ni (t,x)\mapsto I(t,x)$ (see Lemma \ref{P:J0LocalLip}), for
all compact sets $K\subseteq \R_+^*\times\R$, $\sup_{(t,x)\in K} I(t,x)^2<+\infty$.
This completes the whole proof of Lemma \ref{L:InDt}.
\end{proof}

\begin{proof}[Proof of Theorem \ref{T:ExUni}]

The proof follows the same six steps as those in the proof of
\cite[Theorem 2.4]{ChenDalang13Heat} with some minor changes:

Both proofs rely on estimates on the kernel function $\calK(t,x)$. Instead of an explicit formula for the SHE (see
\cite[Proposition 2.2]{ChenDalang13Heat}), Theorem \ref{T:K-bounds} ensures the finiteness of $\calK(t,x)$ and provides a bound on it.

In the Picard iteration scheme, i.e., Steps 1--4 in the proof of \cite[Theorem 2.4]{ChenDalang13Heat},
we need to check the $L^p(\Omega)$-continuity of the stochastic integral,
which then guarantees that at the next step, the integrand is again in $\calP_2$, via \cite[Proposition 3.4]{ChenDalang13Heat}.
Here, the statement of \cite[Proposition 3.4]{ChenDalang13Heat} is still true by replacing in its proof
\cite[Proposition 3.5]{ChenDalang13Heat} by either Proposition \ref{P:G-SD} for the slow diffusion equations or
Proposition \ref{P:G-FD} for the fast diffusion equations, and replacing \cite[Proposition 5.3]{ChenDalang13Heat} by Proposition \ref{P:G-Margin}.

In the first step of the Picard iteration scheme, the following property,
which determines the set of the admissible initial data, needs to be verified:
for all compact sets $K\subseteq \R_+\times\R$,
\[
\sup_{(t,x)\in K}\left(\left[1+J_0^2\right]\star
G_\beta^2 \right) (t,x)<+\infty.
\]
For the SHE, this property is proved in \cite[Lemma 3.9]{ChenDalang13Heat}.
Here, Lemma \ref{L:InDt} gives the desired result with minimal requirements on the initial data.
This property, together with the calculation of  the upper bound on $\calK(t,x)$
in Theorem \ref{T:K-bounds}, guarantees that all the $L^p(\Omega)$-moments of $u(t,x)$ are finite.
This property is also used to establish uniform convergence of the Picard iteration scheme, hence $L^p(\Omega)$--continuity of $(t,x)\mapsto I(t,x)$.

The proof of \eqref{E:SecMom-Lower} is identical to that of the corresponding property in \cite[Theorem 2.4]{ChenDalang13Heat}.
This completes the proof of Theorem \ref{T:ExUni}.
\end{proof}

\addcontentsline{toc}{section}{Bibliography}

\def\polhk#1{\setbox0=\hbox{#1}{\ooalign{\hidewidth
  \lower1.5ex\hbox{`}\hidewidth\crcr\unhbox0}}} \def\cprime{$'$}
  \def\cprime{$'$}

%

\end{document}